\documentclass{amsart}
\newcommand{\bC}{\mathbb{C}}
\newcommand{\bP}{\mathbb{P}}
\newcommand{\bR}{\mathbb{R}}
\newcommand{\bZ}{\mathbb{Z}}

\newcommand{\bT}{\mathbb{T}}
\newcommand{\Cstar}{{\mathbb{C}}^{*}}
\newcommand{\mL}{\mathcal{L}}
\newcommand{\mK}{\mathcal{K}}

\newcommand{\calN}{\mathcal{N}}
\newcommand{\calD}{\mathcal{D}}

\newcommand{\p}{\partial}

\newcommand{\uX}{\underline{\Bbb X}}

\newcommand{\tLM}{\widetilde{{\mathcal{L}}M}}
\newcommand{\tq}{\tilde{q}}
\newcommand{\tqk}{\tilde{q_k}}

\newcommand{\tLNM}{\widetilde{{{\mathcal{L}}_N}M}}
\newcommand{\tLNCn}{\widetilde{{{\mathcal{L}}_N}{{\bC}^n}}}

\newtheorem{theorem}{Theorem}

\newtheorem{proposition}{Proposition}
\newtheorem{example}{Example}

\newtheorem{definition}{Definition}

\newtheorem{conjecture}{Conjecture}

\begin{document}

\title[QUANTUM COHOMOLOGY AND MORSE THEORY ON THE LOOP SPACE]{QUANTUM COHOMOLOGY AND MORSE THEORY ON THE LOOP SPACE OF
TORIC VARIETIES}
\author{Yiannis Vlassopoulos}
 \email{yvlassop@math.uoa.gr}
\date{February 21, 2002}
\address{Dept. of Mathematics\\
         University of Athens\\
         Panepistimiopolis, 15784, Athens, Greece}

\begin{abstract}
On a symplectic manifold $M$, the quantum product defines a
complex, one parameter family of flat connections called the
A-model or Dubrovin connections. Let $\hbar$ denote the parameter.
Associated to them is the quantum $\mathcal{D}$ - module
${\mathcal{D}}/I$ over the Heisenberg algebra of first order
differential operators on a complex torus. An element of $I$ gives
a relation in the quantum cohomology of $M$ by taking the limit as
$\hbar\to 0$. Givental \cite{HomGeom}, discovered that there
should be a structure of a $\mathcal{D}$ - module on the (as yet
not rigorously defined) ${S^1}$ equivariant Floer cohomology of
the loop space of $M$ and conjectured that the two modules should
be equal. Based on that, we formulate a conjecture  about how to
compute the quantum cohomology $\mathcal{D}$ - module in terms of
Morse theoretic data for the symplectic action functional. The
conjecture is proven in the case of toric manifolds with
$\int_d{c_1}> 0$ for all nonzero classes $d$ of rational curves in
$M$.

\end{abstract}
\maketitle

\section{Introduction}
In this paper we will study the quantum cohomology and more
generally the quantum $\calD$ - module structure, of symplectic
toric manifolds by relating it to Morse theory of the unperturbed
symplectic action functional on the loop space. In particular we
will use $S^1$ equivariant cohomology of the action functional.

This program was initiated by Givental in \cite{HomGeom} and
provided the inspiration for the methodology applied later in
\cite{EGWI} in the context of Kontshevich's space of stable maps.

In order to describe the main theorem let us briefly recall a few
things about quantum cohomology. We follow mainly Givental
\cite{EGWI} in this introductory exposition.

Let $(M,\omega)$ be a symplectic manifold and choose also a
compatible almost complex structure $J$. A pseudo-holomorphic
curve is a map $f:{{\bP}^1}\to M$ whose derivative is complex
linear \footnote{To avoid any confusion, let us point out that in
this paper by a pseudo-holomorphic curve we will always mean one
that has arithmetic genus 0.}. Let $d\in {H_2}(M,{\bZ})$ be a
homology class. Kontsevich \cite{Kont} invented the correct space
parameterizing pseudo-holomorphic curves in $M$ with $k$ marked
points. It is called the space of stable maps and we'll denote it
by $M_{k,d}$. Let us consider $M_{3,d}$. Its elements are
equivalence classes of 4-tuples, $(f,{x_1}.{x_2},{x_3})$ where $f$
is the map and the $x_i$ 's are the marked points. The 4-tuple
must satisfy the stability condition that it has at most a
discrete group of automorphisms. Two 4-tuples are equivalent if
there is an automorphism of ${\bP}^1$ that takes one to the other.
The space $M_{3,d}$ is at worst an orbifold (Kontshevitch
\cite{Kont}) if $M$ is convex (i.e., ${H^1}({{\bP}^1},{f^*}TM)=0$
for all stable $f$). \footnote{Main examples of convex spaces are
homogeneous spaces.} Moreover it comes equipped with three
evaluation maps ${ev_i}:{M_{3,d}}\mapsto M$ for $i=1,2,3$ given by
${ev_i}(f,{x_1}.{x_2},{x_3})=f({x_i})$. Let $a,b$ be classes in
$H^{2*}(M,\bC)$. Let $(\ , \ )$ denote the intersection pairing.
Let also ${p_1},\dots {p_r}$ be classes in the K\"{a}hler cone
$\mK$ of $M$ which form a basis of ${H^2}(M,\bZ)$. The K\"{a}hler
cone is the cone in ${H^2}(M,\bR)$ which consists of all classes
whose integral over any (pseudo-) holomorphic curve is
non-negative. Finally let ${d_i}=\int_{d}{p_i}$ and
${q_i}=e^{t_i}$ be complex variables (which can be thought of as
coordinates on a complex torus). The quantum product $a*b$ is
defined by the property that
$$(a*b,c)=\sum_{d}{q^d}\int_{M_{3,d}}{{ev_1}^*}(a)\wedge {{ev_2}^*}(b)
\wedge {{ev_3}^*}(c),$$ where $q^d=\prod_{i=1}^{r}{q_i}^{d_i}$ and
the sum is over all homology classes $d$ of pseudo-holomorphic
curves. The number ${(a*b,c)_d}=\int_{M_{3,d}}{{ev_1}^*}(a)\wedge
{{ev_2}^*}(b) \wedge {{ev_3}^*}(c)$ is called a Gromov-Witten
invariant (of the symplectic structure). It should be thought of
geometrically, as counting the number of curves in homology class
$d$,  meeting classes dual to $a$, $b$ and $c$, when the number of
such curves is finite. Otherwise it is $0$. The (small) quantum
cohomology ring of $M$ is the ring
${SQH^*}(M)={H^{2*}}(M,{\bZ})\bigotimes{\bC}[[{q_1},\dots,{q_r}]]$
equipped with the quantum product.

The quantum product is commutative and associative. This last
property is highly nontrivial and makes for many interesting
consequences by itself.  It turns out that the associativity can
be reformulated as the flatness of the following complex one
parameter family of connections
$$\nabla_{\hbar}={\hbar}d-\sum_{i=1}^r\frac{dq_i}{q_i}\wedge {p_i}*$$
acting on elements of ${SQH^*}(M)$, where $\hbar$ denotes the
complex parameter. This is called the Dubrovin or A-model
connection and as we'll see shortly, it is a more fundamental
object than the quantum product. Givental in his remarkable paper
\cite{EGWI} found a formula for flat sections of $\nabla_{\hbar}$.
It's clear that for flat sections of $\nabla_{\hbar}$, quantum
multiplication by $p_i$ is translated to differentiation and
therefore we may expect that relations in the quantum ring may be
translated to differential equations. This was formulated
explicitly by Givental in the following fashion:
 We may associate  to $\nabla_{\hbar}$ a certain
$\mathcal{D}$ - module $\mathcal{D}/I$ over the algebra of
Heisenberg differential operators. This has the property that if
the operator $D({\hbar}{q_i}\frac{\p}{{\p}{q_i}},{q_i},{\hbar})$
is in the ideal $I$, then the relation $D({p_i}*,{q_i},0)=0$ holds
in the quantum cohomology ring ${SQH^*}(M)$. Therefore the
$\mathcal{D}$ - module appears to be the real quantum object while
the quantum ring arises as its ``semi-classical approximation''
when $\hbar \to 0$\ !

To describe the $\mathcal{D}$ - module $\mathcal{D}/I$ we need to
introduce a new ingredient. This is the line bundle $L$ over
$M_{2,d}$ which is the universal cotangent line at the second
marked point, i.e., the line bundle whose fiber over
$[S,({x_1},{x_2}),f]$ is the cotangent line to $S$  at the second
marked point.

Let $c$ denote the first Chern class of $L$.\footnote{Because of
the meaning it acquires in string theory, $c$ is called a
gravitational descendent.} Now choose bases ${T_0},\dots,{T_m}$
and  ${T^0},\dots,{T^m}$ of ${H^{2*}}(M,{\bZ})$ such that
$({T^i},{T_j})=\delta_{i,j}$. We arrange that ${T_0}=1\in
{H^0}(M,{\bZ})$ and ${T_i}={p_i}$ for $i=1,\dots r$ where $p_i$ as
before. Givental's result \cite{EGWI} is the following: Let $G$ be
the $H^{2*}(M,{\bC})$ valued function defined as:
\begin{equation}
G={e^{plnq/{\hbar}}}(1+\sum_{d}{q^d}{{ev_1}_*}(\frac{1}{{\hbar}-c})),
\end{equation}
where $d$ ranges over all non-zero homology classes of
pseudo-holomorphic curves and ${ev_1}:{M_{2,d}}\to M$ is
evaluation at the first marked point. Then the ideal I is
generated by all polynomial differential operator that annihilate
the components of $G$.

The object of this paper is to compute the quantum $\mathcal{D}$ -
module and specifically the function $G$, \textit{not}
\textit{via} \textit{the} \textit{space} \textit{of}
\textit{stable} \textit{maps}, \textit{but} \textit{rather}
\textit{by} \textit{working} \textit{on} \textit{the}
\textit{loop} \textit{space} \textit{of} $M$. This is desirable
for several reasons, one of them being that the set up seems to be
more natural (at least as far as the $\mathcal{D}$ - module is
concerned) and therefore we may get a better geometric
understanding of the formulas . Another one is that we have a
built in $S^1$ symmetry which, if correctly understood, should
simplify the problem. (By definition the Kontsevich space of
stable maps in $M$ doesn't have that symmetry and Givental in
\cite{EGWI} uses instead the space of stable maps $f:{\bP^1} \to M
\times {\bP^1}$ and the $S^1$ action on the second factor.) Both
of these observations are materialized here to some extent and we
expect they will even more in the future. Yet another reason to
insist with the loop space is that since its $S^1$ equivariant
cohomology is computed by (a version of) the cyclic bar complex
\cite{GJP} we expect fruitful interaction and more powerful
calculational tools to emerge when the problem is properly
formulated in that setting.

The connection with the loop space was first explained by Givental
in \cite{HomGeom}. The connection between pseudo-holomorphic
curves and the loop space was already present in the work of Floer
and his celebrated proof of the Arnold conjecture which resulted
in the definition of Floer homology. To explain the appearance of
the loop space let us start with the space $\mathcal{L}M$ of free
contractible loops in $M$. We can define the action functional $H$
by
$$H(\gamma)=\int_{D_{\gamma}}{\omega},$$
where $\gamma$ is a contractible loop and $D_{\gamma}$ a disc
contracting it. It is multi-valued if there are homologically
non-trivial spheres. To resolve the ambiguity we lift it to the
covering space $\tLM$ of $\mathcal{L}M$, with covering group the
group of spherical classes in $M$. Assume for simplicity that $M$
is simply connected, then ${H_2}(M,{\bZ})$ is generated by
spherical classes. Now the key is, that $H$ has the remarkable
property that its flow lines are pseudo-holomorphic cylinders!
Moreover $H$ is a Hamiltonian function with respect to the obvious
circle action on $\tLM$ and the symplectic form induced from the
symplectic form on $M$. The critical manifolds correspond to
trivial loops and are copies of $M$, one for every degree $d\in
{H_2}(M,{\bZ})$, i.e., for every floor of the cover. Denote by
$M_0$ the copy on which $H$ has value $0$ and by $M_d$ its
translation by $d$.

A formal application of the $S^1$ equivariant localization theorem
suggests that the (Floer) $S^1$ equivariant cohomology of $\tLM$
should be simply
${FH^*}_{S^1}(\tLM)={H^*}(M,{\bC}[{\tq},{{\tq}^{-1}}])({\hbar})$
where ${\bC}[{\tq},{{\tq}^{-1}}]$ is the group ring of the
covering group.

Givental's observation is that ${FH^*}_{S^1}(\tLM)$ bears the
structure of a $\mathcal{D}$ - module over the Heisenberg algebra
of differential operators. This is shown by extending the classes
$\{{p_1},\dots ,{p_r}\}$ to equivariant classes $\{{P_1},\dots
,{P_r}\}$ (see (24)). Then if we think of the $P_k$ acting by
multiplication and the ${\tqk}^d$ by pullback it is easy to show
(27) that
$$[{P_j},{\tqk}]={{\delta}_{j,k}}{\hbar}{\tqk}$$
Givental conjectures that this $\mathcal{D}$ - module is the
quantum $\mathcal{D}$ - module.

Now let $\calN_d$ denote the normal bundle to $M_d$. Let also
${\calN_d}^+$ and ${\calN_d}^-$ denote the positive and negative
normal bundles to $M_d$, with respect to $H$. They are both
infinite dimensional. In section four we argue that if we can make
sense of the ratio of $S^1$ equivariant E\"{u}ler classes
$\frac{{e_{S^1}}({{\calN_d}^+})}{{e_{S^1}}({{\calN_0}^+})}$, then
the cohomology valued function
\begin{equation}
F=\sum_{d} {e^{{plnq}/{\hbar}}}{q^d}
\frac{{e_{S^1}}({{\calN_d}^+})}{{e_{S^1}}({{\calN_0}^+})}
\end{equation}
where $d$ ranges over all homology classes of pseudo-holomorphic
curves,
 should generate the quantum $\mathcal{D}$ - module, in the sense
that the ideal $I$ should be generated by operators annihilating
the components of $F$ in some basis. This is formulated as
conjecture (1) in section (4).

Now if $M$ is a toric variety, we construct approximations of
$\tLM$ by finite dimensional spaces of loops of finite but
arbitrarily large modes (in the Fourier expansion). In that case
we are able to prove (Proposition 3, section 5) that the ratios of
E\"{u}ler classes stabize and therefore we may define a
``$\it{stable\ \ ratio\ of\ Euler\ classes}$" (Definition 2,
section 5). Using this stable ratio in the formula for $F$ we are
then able to prove the main theorem (Theorem 4, section 5) which
says that if a $M$ is a positive K\"{a}hler toric manifold in the
sense that $\int_d{c_1}> 0$ for all nonzero classes $d$ of
arithmetic genus $0$ curves in $M$ then
\begin{equation}
F=G
\end{equation}
and therefore \textit{indeed} \textit{generates} \textit{the}
\textit{quantum} $\mathcal{D}$ - \textit{module}. This is the
content of section five.

The structure of the paper is as follows: In the next section we
gather the elements of the theory of quantum cohomology that are
needed and which are contained mainly in \cite{EGWI}. In the third
section we explain the original idea of Givental \cite{HomGeom}
relating the quantum $\calD$ - module and $S^1$ equivariant Floer
homology of the loop space, via a sort of "Fourier" transform of
cycles. Then a heuristic conjecture is presented on how to compute
the "Fourier" transform of cycles arising in Givental's work. In
the fifth section the conjecture is formulated rigorously for the
case of toric manifolds and then proven. Finally in the last
section we make some observations and indicate what seem to be
 interesting directions for extension of this work.

\textit{Acknowledgements}: The results proven in this paper form
part of my Ph.D thesis \cite{thesis}. I would like to express my
heartfelt thanks and gratitude to my advisor David Morrison for
his guidance and insight.

It is a pleasure to thank Mark Stern, Robert Bryant, Richard Hain,
John Harer and Paul Horja for many helpful conversations.  It is
also a pleasure to acknowledge several valuable discussions with
A. Givental.

\section{The Quantum $\calD$ - Module}

We gather here the basic facts about the quantum $\calD$ - module
that will be needed later on. Let $M$ be a K\"{a}hler manifold and
let ${T_0},{T_1},\dots,{T_r},\dots,{T_m}$ be a basis of
${H^{2*}}(M,{\bZ})$ where ${T_0}$ is the identity and
${T_1},\dots,{T_r}$ are classes in the K\"{a}hler cone $\mK$ that
generate ${H^2}(M,\bZ)$. Recall that $\mK$ is the cone of classes
whose integral over holomorphic curves is non-negative. Set
${p_i}={T_i}$ for $i=1,\dots r$.
 The small quantum ring is defined
to be the vector space
$${SQH^*}(M)={H^{2*}}(M,{\bZ})\bigotimes{\bC}[[{q_1},\dots,{q_r}]]$$
equipped with the quantum product $*$ already defined in the
introduction. That is, if $( \ , \ )$ denotes the intersection
pairing and $a,b,c\in {H^{2*}}(M,{\bZ})$ we have
\begin{equation}
(a*b,c)=\sum_{d}{q^d}\int_{M_{3,d}}{{ev_1}^*}(a)\wedge
{{ev_2}^*}(b) \wedge {{ev_3}^*}(c),
\end{equation}
 where
$q^d=\prod_{i=1}^{r}{q_i}^{d_i}$, ${d_i}={\int_d}{p_i}$ and the
sum is over all homology classes $d\in H_2(M,\bZ)$ of holomorphic
curves. $M_{3,d}$ is the Kontsevich space of stable maps with
three marked points, whose image has arithmetic genus $0$ and
degree $d$ in $H_2(M,\bZ)$. $SQH^*(M)$ is graded if we assign
cohomology classes their usual degree  and declare
$deg({{q_1}^{d_1}}\dots{{q_r}^{d_r}})=2\int_{d}{c_1}(TM)$.
 The reason for this grading
is that $\int_{M_{3,d}}{{ev_1}^*}(a)\wedge {{ev_2}^*}(b) \wedge
{{ev_3}^*}(c)$ is $0$ unless the sum of the degrees of $a,b,c$ is
equal to the dimension of $M_{3,d}$ which, as computed in
\cite{Kont}, is: $dim_{\bC}M_{3,d}=dim_{\bC}M +
\int_{d}{c_1}({T_M})$. Recall that the quantum product is
commutative and associative (see \cite{EGWI}).

Introduce now a one parameter family of connections with regular
singular points, depending on the complex parameter $\hbar$ and
defined by
\begin{equation}
\nabla_{\hbar}={\hbar}d-\sum_{i=1}^r\frac{dq_i}{q_i}\wedge {p_i}*.
\end{equation}
$\nabla_{\hbar}$ are called Dubrovin or $A$ model connections,
 and act on power series in the $q_i$ with coefficients in
$H^{2*}(M,{\bZ})$, in other words on elements of $SQH^*(M)$. As
already mentioned, one reason for considering these connections is
that flat sections, if they exist, will provide a passage from
quantum product to differentiation and from relations in the
quantum ring to differential equations. Moreover, motivation for
introducing the connections comes from mirror symmetry. In fact in
the Calabi-Yau case they are the counterpart of the Gauss-Manin
connection corresponding to the mirror family. On the other hand,
the reason for having a whole pencil of connections is best
understood from the point of view of the loop space and will be
explained in the next section.

Now in Givental \cite{EGWI} it is proven that:

\begin{proposition}
The connection $\nabla_{\hbar}$ is flat for any value of $\hbar$.
\end{proposition}

 The fact that $\nabla_{\hbar}$ is
flat means that we can find flat sections. That is, sections $s$
such that $\nabla_{\hbar}s=0$.

One of the remarkable results of  \cite{EGWI} is the explicit
computation of the flat sections. To describe them we introduce
first the line bundle $L$ over $M_{2,d}$ which is the universal
cotangent line at the second marked point, i.e., the line bundle
whose fiber over $[S,({x_1},{x_2}),f]$ is the cotangent line to
$S$  at the second marked point. Let $c$ denote the first Chern
class of $L$. Choose now basis ${T^0},\dots,{T^m}$ of
${H^{2*}}(M,{\bZ})$ such that $({T^i},{T_j})=\delta_{i,j}$. We
still have that ${T_0}=1\in {H^0}(M,{\bZ})$ and that ${p_i}={T_i}$
for $i=1,\dots r$ where ${p_i} \in \mK$ are  chosen to be a basis
of ${H^2}(M,{\bZ})$. Givental's result \cite{EGWI} is:
\begin{theorem}
The sections
$$s_{\beta}={e^{{plnq}/{\hbar}}}{T_{\beta}}+
\sum_{\alpha}{T^{\alpha}}\sum_{d\in {\check{\mathcal{K}}},{d\neq
0}}{q^d}<{T_{\alpha}},{e^{{plnq}/{\hbar}}}\frac{T_{\beta}}{{\hbar-c}}>_d$$
for $\beta=0,\dots m$ are flat and they provide a basis of the
space of flat sections.
\end{theorem}
 Here $q^d$ is notation for
${{q_1}^{d_1}}\dots{{q_r}^{d_r}}$, $plnq$ is notation for
${p_1}{ln{q_1}}+\dots +{p_r}{ln{q_r}}$ and
$$<{T_{\alpha}},{e^{{plnq}/{\hbar}}}\frac{T_{\beta}}{{\hbar-c}}>_d=
\int_{M_{2,d}}{{ev_1}^*}(T_{\alpha})\wedge
\frac{{{ev_2}^*}({e^{{plnq}/{\hbar}}}{T_\beta})}{\hbar-c},$$ where
$M_{2,d}$ is the space of arithmetic genus $0$ and degree $d$
stable maps with two marked points. Finally $\check{\mathcal{K}}$
is the cone in $H_2(M,\bZ)$ consisting of classes of holomorphic,
arithmetic genus $0$ curves. It is dual to the K\"{a}lher cone
$\mathcal{K}$.  Note also that the matrix:
\begin{equation}
{s_{{\alpha},{\beta}}}=({T_{\alpha}},s_{\beta})=({T_{\alpha}},{e^{{plnq}/{\hbar}}}{T_{\beta}})+
\sum_{d\in {\check{\mathcal{K}}},d\ne
0}{q^d}<{T_{\alpha}},{e^{{plnq}/{\hbar}}}\frac{T_{\beta}}{{\hbar-c}}>_d.
\end{equation}
 is the fundamental solution matrix of the flat section equation.

Let us now explain the relation of the $A$-connection to the small
quantum ring. Let $G$ be the following function with values in
$H^{2*}(M,{\bC})$ :
\begin{equation}
G={e^{plnq/{\hbar}}}(1+\sum_{d\in {\check{\mathcal{K}}},d\ne
0}{q^d}{{ev_1}_*}(\frac{1}{{\hbar}-c})),
\end{equation}
where ${ev_1}:{M_{2,d}}\to M$ is evaluation at the first marked
point. Then $G$ has the property (and is determined by it):
$(G,{T_{\beta}})=(1,s_{\beta})$. Indeed,
$$(G,{T_{\beta}})=({e^{{plnq}/{\hbar}}},{T_{\beta}})+
\sum_{d\in {\check{\mathcal{K}}},d\ne
0}{q^d}\int_{M}{e^{{plnq}/{\hbar}}}{{ev_1}_*}(\frac{1}{{\hbar}-c})
\wedge {T_{\beta}}=$$
$$=(1,{e^{{plnq}/{\hbar}}}{T_{\beta}})+\sum_{d\in {\check{\mathcal{K}}},d\ne
0}{q^d}\int_{M_{2,d}}\frac{1}{{\hbar}-c}\wedge
{{ev_1}^*}({e^{{plnq}/{\hbar}}}{T_\beta})=(1,s_{\beta})$$
Therefore we have that
\begin{equation}
G=\sum_{\beta}(1,s_{\beta}){T^{\beta}}.
\end{equation}
Recall that ${T_0}=1$ and therefore the components of $J$ form the
first row of the solution matrix $(s_{{\alpha},{\beta}})$.

The following proposition is due to Givental \cite{EGWI}
\begin{proposition}
Let $D({\hbar}{q_i}\frac{\p}{{\p}{q_i}},{q_i},{\hbar})$ be a
polynomial differential operator that annihilates the components
of $G$. Then the relation $D({p_i}*,{q_i},0)=0$ holds in
$SQH^*(M)$.
\end{proposition}

Let $\mathcal{D}$ denote the Heisenberg algebra of differential
operators on holomorphic function on a torus with coordinates
${q_i}={e^{t_i}}$. It is by definition generated by the operators
${\hbar}{q_i}\frac{\p}{{\p}{q_i}}={\hbar}\frac{\p}{{\p}{t_i}}$ and
multiplication by ${q_i}={e^{t_i}}$. Let $I$  be the ideal of all
polynomial differential operators
$D({\hbar}{q_i}\frac{\p}{{\p}{q_i}},{q_i},{\hbar})$ that
annihilate the components of $G$.
\begin{definition}
The $\mathcal{D}$ - module $\mathcal{D}/I$ is called the {\itshape
quantum cohomology $\mathcal{D}$ - module} of $M$.
\end{definition}

The proposition above shows that the real quantum object is the
$\calD$ - module or equivalently the $A$ model connection, while
the quantum ring should be considered as the semi-classical limit
where $\hbar\to 0$ . Our objective is to compute the $\calD$ -
module in terms of the loop space of $M$. We shall turn to this
next.

\section{Equivariant Floer theory}

 Lets start first by considering the $S^1$ equivariant Floer
homology of the unperturbed action functional $H$ in the case of a
general symplectic manifold which is not necessarily toric.

Let $(M,\omega)$ be a compact symplectic manifold. Let $J$ be a
compatible or calibrated almost structure on $M$. By this we mean
that $\omega(v,Jv)\geq 0$ for all nonzero $v\in TM$ and
$\omega(Jv,Jw) = \omega(v,w)$. The symplectic form $\omega$ along
with $J$ define an invariant metric $g$ on $TM$ by
$g(v,w)=\omega(v,Jw)$.
 Let
$\mathcal{L}M$ be the space of smooth maps $\gamma:S^1 \to M$ such
that $\gamma(S^1)$ is contractible. We call $\mathcal{L}M$ the
{\itshape loop space} of $M$. The loop space inherits a symplectic
structure $\Omega$ and an almost complex structure which we shall
denote also by $J$. To describe them lets first consider the
tangent bundle $T\mL M$. The tangent space of $\mL M$ at a loop
$\gamma$ is ${T_{\gamma}}\mL M=\mathcal{C}({\gamma}^*TM)$, where
$\mathcal{C}$ denotes the space of sections. In other words an
element of ${T_{\gamma}}\mL M$ is a vector field along the loop
$\gamma$. Consider now the K\"{a}hler cone ${\mK} \subset
{{H^2}(M,{\bR})}$ of $M$. We have defined $\mK$ to be the cone of
classes in ${H^2}(M,{\bR})$ whose integral over any
pseudo-holomorphic curve is greater than or equal to zero. Assume
that  $\mK$ is spanned by the classes of symplectic two forms
 $\omega_1,\dots,\omega_l$.
Let $v$ and $w$ be elements of $T_{\gamma}{\mL}M$  then we define:
\begin{equation}
{{\Omega_k}|_{\gamma}}(v,w)=\int_{S^1}{\omega_k}(v(t),w(t))dt.
\end{equation}
It is not hard to show that the ${\Omega_k}$ are also symplectic.
Moreover $J$ induces an almost complex structure by
$(Jv)(t)=J(v(t))$. Finally $T_{\gamma}{\mL}M$ becomes pre-Hilbert
with the inner product
\begin{equation}
g_{\gamma}(v,w)={\Omega}_{\gamma}(v,Jw),
\end{equation}
where
$${{\Omega}_{\gamma}}(v,w)=\int_{S^1}{\omega}(v(t),w(t))dt.$$
Introduce now {\itshape action functionals}
\begin{equation}
{H_k}(\gamma)=\int_{D_{\gamma}}{{\omega}_k},
\end{equation}
for $k=1\dots l$ and
\begin{equation}
{H}(\gamma)=\int_{D_{\gamma}}{\omega},
\end{equation}
where $D_{\gamma}$ is a disk contracting the loop $\gamma$. These
are in general not well defined since different disks contracting
the same loop will not have the same symplectic areas. The
ambiguity in $H_k$ is clearly given by the periods
$$\int_{S}{\omega_k},$$
where $S$ is a sphere obtained by gluing two different disks
contracting $\gamma$, along their common boundary. The functions
$H_k$ become well defined only on the covering of $\mL M$ with
group of deck transformations the group of spherical periods of
the symplectic forms ${{\omega}_1},\dots ,{{\omega}_l}$. We shall
denote this space by $\tLM$. We can describe $\tLM$ explicitly as
equivalence classes of pairs $({\gamma},g)$ where $\gamma:{S^1}\to
M$ is a loop and $g:D\to M$ is such that ${g|_{\p D}}={\gamma}$.
Define $({\gamma},{g_1}) \sim ({\gamma},{g_1})$ if and only if
$g_1 \# (-g_2)$ represents a class $A \in H_2(M,\bZ)$ such that
$\int_{A}{\omega_k}=0$ for all $k=1\dots l$. Observe that by
definition $\tLM$ carries an action, denoted by $\cdot$ of the
group $\Gamma$ of spherical classes in ${H_2}(M,{\bZ})$ such that
$(A\cdot ({\gamma},g))\#(\gamma,-g)=-A$ for all $A\in {\Gamma}$.
Notice that we could have chosen a positive instead of a negative
sign in the definition of the action of ${\Gamma}$. The reason for
our choice will become apparent later (see footnote (8), page 13).
Note also that since $\omega$ is a linear combination of
${{\omega}_1},\dots ,{{\omega}_l}$ it follows that $H$ is the same
linear combination of ${H_1},\dots {H_k}$ and therefore also
becomes a well defined function on $\tLM$. Now it is not hard to
compute that
\begin{equation}
dH_k|_{\gamma}(v)=-\int_{S^1}{{\omega}_k}(\dot{\gamma},v(t))dt,
\end{equation}
where $\dot{\gamma}$ denotes the vector field tangent to $\gamma$.
To see this, let $v(t)$ be a vector field along $\gamma$ and let
${{\gamma}_s}({e^{it}})$, $s\in \bR$ be a curve representing $v$.
This means that ${{\gamma}_0}=\gamma$ and
$\frac{d}{ds}|_{s=0}{{\gamma}_s}=v(t)$. Moreover we can arrange
that ${\gamma}_s({e^{it}})=a$ for all $s\leq -1$ where $a$ is a
point in $M$. Now consider the family of maps on the disc $D$,
contracting the family of loops, given by
${u_s}({e^{r+it}})={{\gamma}_{s+r}}({e^{it}})$ for $r\leq 0$.
Notice that ${u_s}(e^{it})={{\gamma}_s}({e^{it}})$ so that $u_s$
maps the boundary of $D$ to the image of ${\gamma}_s$. Moreover
${u_0}({e^{r+it}})={{\gamma}_{r}}({e^{it}})$. Finally ${u_s}=a$
for $s\leq -1$. Now
$${dH_k}|_{\gamma}(v)={\frac{\p}{\p s}|_{s=0}}{H_k}({\gamma}_s)=
{\frac{\p}{\p s}|_{s=0}}\int_{D}{{u_s}^*}{{\omega}_k}=$$
$$\frac{\p}{\p s}|_{s=0}\int_{D}
{{\omega}_k}(\frac{{\p}{u_s}}{{\p}r},\frac{{\p}{u_s}}{{\p}t})dr\wedge
dt.$$
Finally, applying Stokes' theorem gives (13).

Notice further that $\mL M$ and therefore $\tLM$, support an
obvious $S^1$ action $({e^{i\theta}},{\gamma}({e^{i\phi}}))\mapsto
{\gamma}(e^{i(\theta+ \phi)})$. If we let $X$ denote the vector
generating the Lie algebra of $S^1$ and $\uX$ the induced vector
field on $\tLM$ then we have
\begin{equation}
\uX({\gamma})=\dot{\gamma}.
\end{equation}
Equations (9),(13) and (14) reveal the remarkable fact that
\begin{equation}
i_{\uX} {{\Omega}_k}=-dH_k,
\end{equation}
for $k=1\dots l$  and
\begin{equation}
i_{\uX} {\Omega}=-dH.
\end{equation}
In that case $H_k$ is called a Hamiltonian function for $\Omega_k$
and $H$ a Hamiltonian for $\Omega$. $\uX$ can be thought of as a
symplectic gradient of $H$.

Consider now the flow of $H$. Let $u(s,t):{\bR}\times {S^1}\to M$
be a flow line. Specifically this means that
\begin{equation}
\frac{\p u}{\p s}=\nabla H_{{u_s}(t)},
\end{equation}
where ${u_s}(t)$ is simply $u(s,t)$. On the other hand since
$\Omega$ and $J$ are compatible and the metric on $\tLM$ is given
by (10) we have that
\begin{equation}
\nabla H=-J{\uX}.
\end{equation}
Equations (11) and (12) imply then
\begin{equation}
\frac{\p u}{\p s}=-J{\uX}({{u_s}(t)})=-J\frac{\p u}{\p t}.
\end{equation}
Therefore
\begin{equation}
\frac{\p u}{\p s}+J\frac{\p u}{\p t}=0.
\end{equation}
Now recall that a map between almost complex manifolds is called
(pseudo-) holomorphic if its differential respects the almost
complex structures. It is easy to see that equation (20) means
exactly that $u$ is holomorphic with respect to $J$ on $M$ and the
standard complex structure $j$ on the infinite cylinder
${\bR}\times{S^1}$. To verify that, recall that $j$ is defined by
$$j(\frac{\p}{\p s})=\frac{\p}{\p t} \ and
 \ j(\frac{\p}{\p t})=-\frac{\p}{\p s}.$$
Now $u$ is holomorphic if $du\circ j=J\circ du$. Evaluating at
$\frac{\p}{\p t}$ gives
$$duj(\frac{\p}{\p t})=Jdu(\frac{\p}{\p t}) \iff$$
$$-\frac{\p u}{\p s}=J\frac{\p u}{\p t}.$$

So indeed $u(s,t):{\bR}\times {S^1}\to M$ is a (pseudo-)
holomorphic cylinder!\footnote{In fact we also see that the flow
of $-H$ gives antiholomorphic cylinders.} This is the key reason
why quantum cohomology is related to the loop space.

Floer theory is Morse theory for the action functional $H$ on
$\tLM$. Notice that the critical manifolds are copies of $M$, one
of them corresponding to trivial loops and the rest translations
by the action of the group of deck transformations, i.e., the
group of spherical classes in $H_2(M,\bZ)$. This is easy to see
using for example (16) which identifies the critical manifolds as
the fixed manifolds of the circle action.

 Now the fact that $\tLM$
is infinite dimensional pauses several hard problems one needs to
overcome in order to get a well defined theory. For example, for
any critical manifold both the negative and positive normal
bundles of $H$ are infinite dimensional. Therefore the usual
notion of index doesn't make sense. Moreover the standard Morse
theoretic method of analyzing the topology of a space simply
doesn't work. This is because we cannot  describe the change in
topology when going through a critical manifold, by a gluing of
the negative normal sphere bundle since this is trivial!
\footnote{ The negative normal bundle is trivial since any
infinite dimensional bundle is trivial. Moreover the infinite
dimensional sphere is homotopicaly trivial.} It was Floer's idea
to overcome this problem by constructing a Witten type Morse
theory where the index is defined by counting orbits connecting
critical manifolds. The key point to doing this in this case, is
to use orbits of {\itshape bounded energy}. In other words if
$u(s,t):{\bR}\times {S^1}:\to M$ is a flow line, ie satisfies
(20), then define the energy of $u$ by :
\begin{equation}
E(u)=\frac{1}{2}\int_{0}^{1}\int_{-\infty}^{\infty} {|\frac{\p
u}{\p s}|^2}+{|\frac{\p u}{\p t}|^2}.
\end{equation}
We say that $u$ has bounded energy if $E(u)$ is finite. In fact it
is easy to compute that when $u$ is (pseudo-)holomorphic as is the
case for flow lines, then
$$E(u)=\int_{0}^{1}\int_{-\infty}^{\infty}{u^*}{\omega}.$$
Still to get a well behaved theory we have to perturb the flow
equation by an extra term using a periodic Hamiltonian. The
critical manifolds then become points and the theory can be used
to prove the well known Arnold conjecture for periodic
Hamiltonians.

Floer \cite{Floer} was able to rigorously construct a homology
theory, now called Floer homology,\footnote{Floer constructed the
theory and proved the Arnold conjecture for so called monotone
manifolds (this means that the first Chern class of the manifold
is a positive multiple of the symplectic form). For general
symplectic manifolds the theory was constructed in
\cite{Liu-Tian}.} using these perturbed holomorphic cylinders
connecting periodic orbits. He then showed that Floer homology is
isomorphic, with respect to additive structure, to the singular
homology of $M$ with coefficients in an appropriate ring of
Laurent series. It should be mentioned at this point that in
\cite{RT}, \cite{LT} and \cite{PSS} it is proved (in each work
with different methods) that in fact the isomorphism between Floer
and singular homology, respects the ring structures if
${H^{2*}}(M,\bZ)$ is equipped with the {\itshape Quantum product},
and Floer cohomology with the so called, pair of pants product.

 The unperturbed Morse-Bott-Floer theory
has been worked out to a certain extent by Ruan and Tian in
\cite{RT}.\footnote{In fact in \cite{RT} the authors study a more
general situation where the action functional is perturbed but
without necessarily insisting that the critical manifolds be
points}

Following their paper, the space of connecting orbits between two
critical levels that differ by a class $d\in H_2(M,\bZ)$ should be
taken to consist of maps
$$u:{\bR}\times {S^1}\to M$$
that

1.\ are $J$- holomorphic

2.$$E(u)=\int_{{\bR}\times {S^1}}{u^*}{\omega}<\infty$$

3.$$\lim_{s \to -{\infty}}u(s,t)=point$$

and

$$\lim_{s \to
  \infty}u(s,t)=point,$$

in other words the infinite cylinder closes up at the ends to give
a sphere with two point removed.

4.The homology class of the image of $u$ is $d$.

If we denote the set of such maps $u$ by ${\mathcal{M}}_d$ then we
expect this space to have the same dimension as the space of
holomorphic spheres of degree $d$, the expected dimension of
which, is
$$dim\ M +\int_{d}{c_1}(TM).$$
Indeed the calculation of Morrison \cite{Dave}(p. 277) shows that
if ${\phi}:{S^2}\to M$ is holomorphic and M is a complex manifold
then
\begin{equation}
{\chi}({{\phi}^*}{TM})=
h^0({S^2},{{\phi}^*}{TM})-h^1({S^2},{{\phi}^*}{TM})= dim\ M
+\int_{d}{c_1}(TM).
\end{equation}

In any case  we wish to consider Floer $S^1$ \textit{equivariant}
\textit{cohomology} of $\tLM$ so following Givental
\cite{HomGeom}, we bypass all that and try to use localization
technics instead. The localization theorem relates the equivariant
cohomology of a space with a torus acting on it, to that of the
fixed components of the action. One way \cite{AB} of proving this
theorem rests on an analysis of the
${H^*}(point)={H^*}(\bP^\infty)=\bC[\hbar]$ - module structure of
the equivariant cohomology ring. We refer to this paper or
\cite{Audin} for the more general statement which refers to a
torus action. For our purposes we only need the $S^1$ case. The
result \cite{AB} then is that:
\begin{theorem}
Let $X$ be an $S^1$ (finite dimensional) compact manifold. Let $F$
denote the (possibly disconnected) fixed manifold of the action
and let $i:F\to X$ be the inclusion map. Then
$${i^*}:{{H_{S_1}}^*}(M)\to {{H_{S^1}}^*}(F)$$
induces  an (additive) isomorphism after localization to the field of
rational functions ${\bC}({\hbar})$.
\end{theorem}

Notice that since $F$ is fixed it follows that
${{H_{S^1}}^*}(F)={H^*}(F)\otimes{\bC[\hbar]}$. So the meaning of
this theorem is that, if the fixed manifolds are $\{F_{\alpha}\}$
then there is an isomorphism
$$\Phi:{{H_{S_1}}^*}(M)\to
{\bigoplus_{\alpha}}{H^*}({F_{\alpha}},{\bC}({\hbar}))$$
and so $$\Phi(a)=\sum_{\alpha}{{\lambda}_{\alpha}}{C_{\alpha}}$$
where ${C_{\alpha}}\in {H^*}({F_{\alpha}})$ and ${\lambda}_{\alpha}$
is a {\itshape rational} function in $\hbar$

The more precise statement described in \cite{AB} shows that we
don't really need to invert $\hbar$ to get an isomorphism but
inverting some multiple of it (determined as can be expected by
the equivariant E\"{u}ler class of the normal bundle of $F$) is
enough. In the general torus case the theorem describes precisely
the localization needed in order to get an isomorphism. We refer
to \cite{AB} for the proof.

We would like now to apply this to the space $\tLM$. Since it is
infinite dimensional this is only a {\itshape formal} application, not
a rigorous one.

With this qualification, since the fixed manifolds are copies of
$M$, we expect, after Givental \cite{HomGeom}, the $S^1$
equivariant Floer cohomology of $\tLM$ to be as an additive
object,
\begin{equation}
{FH^*}_{S^1}({\tLM})={H^*}(M,{\bC}[{\tq},{{\tq}^{-1}}]({\hbar})),
\end{equation}
where ${\bC}[\tq,{\tq^{-1}}]$ is notation for the group ring
$\Lambda$ of the group $\Gamma$ of spherical classes in
$H_2(M,\bZ)$. In other words, instead of using a direct sum
notation, we have used the group ring to enumerate the fixed
components of the action.

To be more specific, elements of the ring are formal series
$${\lambda}=\sum_{d\in {\Gamma}}{{\lambda}_d}{{\tq}^d}$$
where ${{\lambda}_d}\in {\bC}$, ${{\tq}^d}=e^{2{\pi} id}$ and we
declare $deg\ {{\tq}^d}=2\int_{d}{c_1}(TM)$.

Assume $M$ is K\"{a}hler and simply connected. $M$ being simply
connected implies that ${\Gamma}={H_2}(M,\bZ)$. Choose a basis
$\{{p_1},\dots ,{p_r}\}$ of the K\"{a}hler cone and let
$\{{A_1},\dots ,{A_r}\}$ be the dual basis of
${\Gamma}={H_2}(M,{\bZ})$ in the sense that
$\int_{A_j}{p_k}=\delta_{j,k}$. Then if $d=\sum_{i=1}^r{d_i}{A_i}$
we let ${{\tqk}}={e^{2{\pi}i {A_k}}}$ and
${{\tq}^d}={e^{2{\pi}id}}=\prod_{k=1}^{r}{{\tqk}^{d_k}}$. In this
fashion the group ring $\lambda$ can be identified with the ring
of formal Laurent series
${\Lambda}={\bC}[[{{\tilde{q}}_1},{{{\tilde{q}}_1}^{-1}},\dots
{{\tilde{q}}_r},{{{\tilde{q}}_r}^{-1}}]]$.

Now let $p_k$ be represented by an $S^1$ invariant k\"{a}hler form
${\omega}_k$. Recall we have defined associated k\"{a}hler classes
${\Omega}_k$  on ${\mathcal{L}}M$. Denote by the same name the
pullbacks on $\tLM$. Let ${d_{\hbar}}=d+{\hbar}{i_{\uX}}$ be the
Cartan differential, where $\uX$ is defined by (14). Introduce now
the equivariant differential forms
\begin{equation}
{P_k}={{\Omega}_k}+{\hbar}{H_k},
\end{equation}
and
\begin{equation}
P={\Omega}+{\hbar}{H}.
\end{equation}
Then equations (15) and (16) imply that
$${d_{\hbar}}{P_k}={d_{\hbar}}{P}=0$$
for $k=1\dots r$.

By definition of the $A_j$ we have:
$$\int_{A_j}{{\omega}_k}={\delta}_{j,k}.$$
Recall we have defined
$${H_k}({\gamma},g)=\int_{D}{g^*}{{\omega}_k}.$$
Recall also that $\Gamma$ acts on $\tLM$ as the group of covering
transformations (so in fact ${\tLM}/{\Gamma}={\mathcal{L}}M$). Now
if we identify ${\tq}^d$ with the covering transformation
corresponding to $d\in {\Gamma}$ we have that
\begin{equation}
{{{\tilde{q}}_j}^*}{H_k}({\gamma},g)=
{H_k}({{\tilde{q}}_j}\cdot ({\gamma},g))={H_k}(({\gamma},g))-
\int_{A_j}{{\omega}_k}={H_k}(({\gamma},g))-{{\delta}_{j,k}},
\end{equation}
where ${{{\tilde{q}}_j}^*}$ denotes the pullback.

Moreover, if we denote by $P_k$ wedge product by the equivariantly
closed form $P_k$ and also denote simply by $\tqk$ the action of
$\tqk$ by pullback then we claim that
\begin{equation}
[{P_j},{\tqk}]={{\delta}_{j,k}}{\hbar}{\tqk}.
\end{equation}
The proof is a simple calculation.
 First notice that:
$${\tqk}(P_j)={{\Omega}_j}+{\hbar}{{{\tilde{q}}_k}^*}{H_j}=
{{\Omega}_j}+{\hbar}{H_j}-{\hbar}{\delta}_{j,k}{\hbar}={P_j}-
{{\delta}_{j,k}}{\hbar}.$$
Next let $\alpha$ be an equivariant
form then
$$[{P_j},{\tqk}]{\alpha}={P_j}{\tqk}{\alpha}-{\tqk}(P_j){\wedge}{\tqk}{\alpha}=
{P_j}{\tqk}{\alpha}-({P_j}-{\delta}_{j,k}{\hbar}){\wedge}{\tqk}{\alpha}=
{\delta}_{j,k}{\hbar}{\tqk}{\alpha}.$$
Besides (27) the operators
$P_j$ and $\tqk$ for $j,k=1\dots r$ satisfy the relations
\begin{equation}
[{P_j},{P_k}]=[{{\tilde{q}}_j},{\tqk}]=0.
\end{equation}

Now let ${t_1},\dots ,{t_r}$ be coordinates on ${\bC}^r$. Then we
have
\begin{equation}
[{\hbar}\frac{\p}{{\p}{t_j}},{e^{t_k}}]={{\delta}_{j,k}}{e^{t_k}},
\end{equation}
where  ${e^{t_k}}$ is thought of as an operator acting by
multiplication on functions of ${e^{t_1}},\dots {e^{t_r}}$ and
${\hbar}\frac{\p}{{\p}{t_j}}$ for $j=1\dots r$  also act on such
functions. The algebra $\mathcal{D}$ of operators generated by
${e^{t_1}},\dots {e^{t_r}}$ and ${\hbar}\frac{\p}{{\p}{t_1}},\dots
,{\hbar}\frac{\p}{{\p}{t_r}}$  is called the Heisenberg algebra of
differential operators. Relations (27) and (28) say that the $S^1$
equivariant Floer cohomology ${{FH_{S^1}}^*}({\tLM})$ carries the
structure of a module \footnote{The reason for the choice of sign
in the action of $\Gamma$ on $\tLM$ is precisely so that we end up
with $e^{t_j}$ instead of $e^{-{t_j}}$.} over the Heisenberg
algebra $\mathcal{D}$ !

In our discussion of the $A$ model connection, in the previous
section, we also encountered a $\mathcal{D}$ - module. That one
consisted of operators which kill the first row of the solution
matrix of the flat section equation for the $A$ connection.
Givental's conjecture is that {\itshape the two $\mathcal{D}$ -
modules are in fact the same}!

Of course there is no chance of proving this unless a rigorous
$S^1$ equivariant Floer theory of the unperturbed action
functional is constructed. In case $M$ is a toric variety though,
we will construct a model for the space $\tLM$ in section (5). If
$M$ is also positive, then we shall be able to prove that the
$\calD$ - module of our model is indeed the same as the quantum
$\calD$ - module.

Notice than in the previous section, we used coordinates ${q_j}$
which are related to the $t_j$ by ${q_j}=e^{t_j}$. It is also
clear, that series in the $q_j$ with values in ${H^{2*}}(M,\bC)$
can be thought of as sections of a trivial bundle with fiber
${H^{2*}}(M,\bC)$, over the (algebraic) torus obtained by the
lattice ${H^2}(M,{\bZ})$ (via complexification and
exponentiation). In particular this torus can be thought of as the
(affine) toric variety associated to a fan consisting of a single
cone, namely the K\"{a}hler cone of $M$. The $q_j$ are then
identified with the toric coordinates.

Now having a $\mathcal{D}$ - module how can we associate a flat
connection ?

Recall that from (23) we have:
\begin{equation}
{FH^*}_{S^1}({\tLM})=
{H^*}(M,{\bC}[{{\tilde{q}}_1},{{{\tilde{q}}_1}^{-1}}\dots ,
{{\tilde{q}}_r}{{{\tilde{q}}_r}^{-1}}]({\hbar})).
\end{equation}

Therefore ${H^{2*}}(M,\bZ)$ is embedded in ${FH^*}_{S^1}({\tLM})$.
Consider again the basis ${T_0},\dots ,{T_m}$ of
${H^{2*}}(M,\bZ)$. We then have
\begin{equation}
{P_k}\wedge ({T_0},\dots ,{T_m})={{\tilde{A}}^k}({T_0},\dots ,{T_m}),
\end{equation}
where ${{\tilde{A}}^k}$ is a matrix with coefficients functions of
$\tqk$ and $\hbar$. The coefficient functions are expected to be
holomorphic so they will not contain any of the ${\tqk}^{-1}$.
Define now a pencil of connections $\tilde{\nabla}_{\hbar}$ acting
on series in the $\tqk$ with values in ${H^{2*}}(M,\bC)$ by:
\begin{equation}
{\tilde{{\nabla}}_{\hbar}}={\hbar}d-\sum_{k=1}^{r}\frac{d{\tqk}}{\tqk}{P_k}{\wedge}.
\end{equation}

The connection $\tilde{\nabla}_{\hbar}$ is expected to be equal to
the $A$ model connection ${\nabla}_{\hbar}$ considered in our
discussion of quantum cohomology based on stable maps. This means
for example that $\tilde{\nabla}_{\hbar}$ should be flat, i.e.,
that flat sections should exist. If
$\sigma=\sum_{j=0}^{m}{f_j}{T_j}$ is a section then
$${\tilde{\nabla}_{\hbar}}{\sigma}=0$$
is equivalent to the system
\begin{equation}
{\hbar}{\tqk}\frac{\p}{{\p}{\tqk}}{{({f_0},\dots ,{f_m})}^t}=
{{\tilde{A}}^k}{{({f_0},\dots ,{f_m})}^t}\ \ for\ k=1\dots r.
\end{equation}
Flatness of course means that (33) is integrable. To shed some
more light we note that (33) can equally be written as
\begin{equation}
{P_k}{\wedge}({T_0},\dots ,{T_m}){{({f_0},\dots ,{f_m})}^t}=
({T_0},\dots ,{T_m})
{\hbar}{\tqk}\frac{\p}{{\p}{\tqk}}{{({f_0},\dots ,{f_m})}^t}.
\end{equation}
In other words the $(m+1)$-tuple ${{({f_0},\dots ,{f_m})}^t}$
defines a $\mathcal{D}$ - module homomorphism between
${FH^*}_{S^1}({\tLM})$ and the sheaf $\mathcal{O}$ of holomorphic
functions on the torus. Therefore we can reformulate our
discussion in an invariant fashion by saying that an element in
${Hom_{\mathcal{D}}}({FH^*}_{S^1}({\tLM}),{\mathcal{O}})$ defines
a locally constant sheaf $\mathcal{V}$  over the torus. This sheaf
defines in turn by the standard procedure a flat connection on the
sheaf $\mathcal{U}={\mathcal{V}}\otimes \mathcal{O}$. As an aside
we note that this may remind the reader of the construction of the
Gauss-Manin connection associated to a family of varieties. The
locally constant sheaf is there, the one associated to the
integral cohomology of the fiber. This is no accident since mirror
symmetry identifies, in the case $M$ is Calabi-Yau, the $A$-model
connection with the Gauss-Manin connection of a certain family of
Calabi-Yau manifolds.

Instead of concentrating on the connection lets look now at the
$\mathcal{D}$ - module itself and try to find a presentation or at
least some relations. We have seen up to now that
${FH^*}_{S^1}({\tLM})$ is generated by ${H^{2*}}(M,\bZ)$ over the
ring $\Lambda$. Geometrically an equivariant Floer cycle
associated to an element $T\in {H^{2*}}(M,\bZ) $ can be
constructed as the boundary loops of all holomorphic discs whose
center lies in a cycle representing the Poincar\'e dual of $T$.
Now notice that the standard way to go between the Heisenberg
algebra and its presentation in terms of the ${P_k}$ and $\tqk$ is
via the {\itshape Fourier transform}. This way relations that
involve the later can be transformed to differential equations
that involve the former. Indeed if $\Gamma$ is in
${FH^*}_{S^1}({\tLM})$ and we denote ${t_1}{P_1}+\dots
+{t_r}{P_r}$ by $tP$ then consider the pairing
\begin{equation}
({e^{{tP}/{\hbar}}},{\Gamma})=\int_{\tLM}{e^{{tP}/{\hbar}}}{\Gamma}.
\end{equation}
We claim that
\begin{equation}
({e^{{tP}/{\hbar}}},R(P,{\tq},{\hbar}){\Gamma})=
R({\hbar}\frac{\p}{{\p}{t}},{e^t},{\hbar})
({e^{{tP}/{\hbar}}},{\Gamma}).
\end{equation}
In other words that the map
$${\mathcal{F}}:{FH^*}_{S^1}({\tLM})\to \mathcal{O}$$ given by
$${\mathcal{F}}(\Gamma)=({e^{{tP}/{\hbar}}},{\Gamma})$$
is an element of ${Hom_{\mathcal{D}}}({FH^*}_{S^1}({\tLM}),{\mathcal{O}})$.

Indeed we can do a bit better than that. If ${C_{T}}\in
{FH^*}_{S^1}({\tLM})$ is such that it has the same localization
$T\in {H^{2*}}(M,{\bZ})$ on every critical manifold, then
\begin{equation}
({e^{{tP}/{\hbar}}}{C_T},R(P,{\tq},{\hbar}){\Gamma})=
R({\hbar}\frac{\p}{{\p}{t}},{e^t},{\hbar})
({e^{{tP}/{\hbar}}}{C_T},{\Gamma}),
\end{equation}
and therefore the map
$${{\mathcal{F}}_T}:{FH^*}_{S^1}({\tLM})\to \mathcal{O}$$ given by
\begin{equation}
{{\mathcal{F}}_T}(\Gamma)=({e^{{tP}/{\hbar}}}{C_T},{\Gamma}),
\end{equation}
is an element of ${Hom_{\mathcal{D}}}({FH^*}_{S^1}({\tLM}),{\mathcal{O}})$.

The reason is that
$$({e^{{tP}/{\hbar}}}{C_T},{\tqk}{\Gamma})=
({{\tqk}^{-1}}{e^{{tP}/{\hbar}}}{C_T},{\Gamma})=
{e^{t_k}}({e^{{tP}/{\hbar}}}{C_T},{\Gamma})$$ and

$$({e^{{tP}/{\hbar}}}{C_T},{P_k}{\Gamma})=
({P_k}{e^{{tP}/{\hbar}}}{C_T},{\Gamma})=
({\hbar}\frac{\p}{{\p}{t_k}}{e^{{tP}/{\hbar}}}C,{\Gamma})=
{\hbar}\frac{\p}{{\p}{t_k}}({e^{{tP}/{\hbar}}}{C_T},{\Gamma}).$$

Thus we can find the differential operators and compute solution
by computing $({e^{{tP}/{\hbar}}},{\Gamma})$ if we can write down
$\Gamma$ and compute the integral.

In fact if $\Gamma=\Delta$ where $\Delta$ is the fundamental Floer
cycle corresponding to the fundamental cycle of $M$, i.e., the
cycle of all boundary loops of holomorphic discs in $M$, and the
cohomology of $M$ is generated by classes in ${H^2}(M,{\bZ})$ then
polynomials $R(P,{\tq},{\hbar})$ such that
$R(P,{\tq},{\hbar}){\Delta}=0$ generate all relations. So if $I_0$
is the ideal generated by such polynomials then
${FH^*}_{S^1}({\tLM})={\bC}[P,{\tq},{\hbar}]/{I_0}$. The reason
for this is that if $f({p_1},\dots ,{p_r})$ is a polynomial in the
generators $\{{p_1},\dots ,{p_r}\}$ of ${H^2}(M,{\bZ})$ then the
corresponding Floer cycle is ${{\Delta}_f}=f({P_1},\dots
,{P_r}){\Delta}$ since out of all loops (boundaries of holomorphic
discs) that have their center in $M$, this pics the ones that are
in the cycle Poincar{\'{e}} dual to $f({p_1},\dots ,{p_r})$. It is
clear now that any polynomial ${R_1}(P,{\tq},{\hbar})$ such that
${R_1}(P,{\tq},{\hbar}){\Delta_f}=0$ induces a relation
${R_1}(P,{\tq},{\hbar})f({P_1},\dots ,{P_r}){\Delta}=0$. Therefore
relations stemming form $\Delta$ generate all relations.

Up to this point our discussion of $S^1$ equivariant Floer theory
of the unperturbed action functional has followed Givental's paper
\cite{HomGeom}. We would like now to propose a conjecture about
how to regularize the integral in (38). In the last section we
shall prove a version of it for toric manifolds.

\section{ A conjecture on the regularization of the Fourier transform
of the Floer fundamental cycle}

Recall first that we have chosen a basis $\{{T_0},\dots ,{T_m}\}$
of ${H^{2*}}(M,{\bR})$. We arrange that $T_0=1$. Choose also a
dual basis $\{{T^0},\dots ,{T^m}\}$ of ${H^{2*}}(M,{\bZ})$  such
that $({T_i},{T^j})=\delta_{i,j}$, where the pairing is the
Poincar\'e pairing. Now to compute the integral in (38) for
$T={T_{\beta}}$ and since the integrand is an equivariantly closed
form, we could attempt to formally use a localization theorem in
equivariant cohomology. The theorem we need is a stronger version
of theorem (2) mentioned before and it is due independently to
Berline-Vergne \cite{BV} and Atiyah-Bott \cite{AB}.
\begin {theorem}
Let $\bT$ be a torus acting on a (finite dimensional) compact
manifold $M$ and let ${\alpha}$ be an equivariantly closed form in
the Cartan model. Then
$$\int_M{\alpha}=
\sum_{F} \int_{F}\frac{{\alpha}_{|F}}{{e_{\bT}}(\calN_F)}$$ Where
the sum is over all the fixed components $F$ of the action and
${e_{\bT}}({\calN_F})$ indicates the $\bT$ equivariant E\"{u}ler
class of the normal bundle, $\calN_F$  to the fixed component $F$.
By ${\alpha}_{|F}$ we denote the pullback of $\alpha$ to $F$ by
the inclusion of $F$ into $M$.
\end{theorem}
 Since $\tLM$ is infinite dimensional an application of
this theorem in our case can only be done in a formal fashion.
This formal application gives:
\begin{equation}
{{\mathcal{F}}_{T_{\beta}}}({\Delta})=\int_{\tLM}{e^{{tP}/{\hbar}}}{C_{T_{\beta}}}{\Delta}=
\sum_{d\in {H_2}(M,{\bZ})}
\int_{M_d}{T_{\beta}}{e^{\sum_{k=1}^{r}{t_k}({{{\omega}_k}/{\hbar}}+
{\int_{d}{{\omega}_k}})}}\frac{{\Delta}_{|{M_d}}}{{e_{S^1}}({\calN_d})},
\end{equation}
where we have used the following notation. First recall that the
action functional $H$ is a function on $\tLM$ whose critical
manifolds are the fixed components of the $S^1$  action and
therefore are just copies of $M$. Denote the copy of $M$ such that
${H_{|M}}=0$ by $M_0$. The action of ${\tq}^d$ maps $M_0$ to
another copy of $M$ which we denote by $M_d$. $\calN_d$ denotes
the normal bundle to $M_d$. This is of course an infinite
dimensional bundle. Notice now that $\calN_d$ carries a
representation of $S^1$ (as a sub-bundle of $T{\tLM}_{|M_d}$)  and
splits to the direct sum of line bundles according to the weights
of this representation. The E\"{u}ler class
${{e_{S^1}}({\calN_d})}$ is therefore some infinite product which
in general, will be divergent. \footnote{It is quite interesting
that Atiyah in \cite{Ati} shows that the inverse of the E\"{u}ler
class of the normal bundle to $M$ sitting as the space of trivial
loops inside the loop space ${\mL}M$, can be normalized by $\zeta$
function regularization and turns out to be equal to
${\Hat{A}}(M)$. Another version of this computation can be found
in Jones and Petrack \cite{JP}.} Moreover recall that
${P_k}={{\Omega}_k}+{\hbar}{H_k}$ and
${{{\Omega}_k}_{|{M_d}}}={{\omega}_k}$ and finally
${{H_k}_{|{M_d}}}=\int_{d}{{\omega}_k}$.

Now we would like to understand better the equivariant Floer
fundamental cycle  $\Delta$. Geometrically it is supposed to be
the Poincar\'e dual of the cycle of loops which are boundary
values of holomorphic discs in $M$. We are interested in the
restriction of $\Delta$ to $M_d$. We have already noticed (20)
that flowlines of the action functional are (pseudo-) holomorphic
cylinders. Therefore flowlines departing from $M_d$ are precisely
(pseudo-) holomorphic discs in $M$. It follows that geometrically
the fundamental cycle should be represented as a Morse-Witten
cycle of $H$, by the formal sum of the unstable
manifolds\footnote{The unstable manifold corresponding to a
critical manifold is the manifold of points on the flow lines
departing from that critical manifold.} of the fixed components
$M_d$. In finite dimensional Morse-Bott-Witten theory an unstable
manifold is fibred over the corresponding critical manifold by the
obvious flow map. Moreover a neighborhood of the zero section in
the positive normal bundle over a critical manifold is
diffeomorphic to a neighborhood of the critical manifold in the
unstable manifold\footnote{See for example  appendix 3 of
\cite{AusBra} for a proof of this in the finite dimensional
case.}. Therefore the cycle defined by the restriction of the
class of the unstable manifold, to the critical manifold should be
the E\"{u}ler class of the positive normal bundle. This implies
that \textit{the} \textit{restriction} \textit{of} $\Delta$
\textit{to} $M_d$ \textit{should} \textit{be} \textit{the}
\textit{equivariant} \textit{E\"{u}ler} \textit{class}
${e_{S^1}}({{\calN_d}^+})$ \textit{of} \textit{the}
\textit{positive} \textit{normal} \textit{bundle} ${\calN_d}^+$ to
$M_d$.

Moreover the restriction of $\Delta$ to $M_d$ for $d$ that cannot
be represented by a (pseudo-) holomorphic curve should be zero. In
other words the non-zero contributions come only from
$d{\in}\check{\mathcal{K}}$.\footnote{Classes not in
$\check{\mathcal{K}}$ have negative $\int_{d}{{\omega}_k}$ for
some $k$. This is not possible for classes of (pseudo-)
holomorphic curves since from equation (21) it follows that this
integral is equal (up to a factor of $\frac{1}{2}$) to the energy
of the curve.}

This informal analysis suggests that (39) becomes:
\begin{equation}
{{\mathcal{F}}_{T_{\beta}}}({\Delta})=\int_{\tLM}{e^{{tP}/{\hbar}}}{C_{T_{\beta}}}{\Delta}=
\sum_{d\in {\check{\mathcal{K}}}}
\int_{M_d}{T_{\beta}}{e^{\sum_{k=1}^{r}{t_k}({{{\omega}_k}/{\hbar}}+
{\int_{d}{{\omega}_k}})}}\frac{{e_{S^1}}({{\calN_d}^+})}{{e_{S^1}}({\calN_d})},
\end{equation}
 Moreover, we have
$${e_{S^1}}(\calN_d)={e_{S^1}}(\calN_0)={e_{S^1}}({{\calN_0}^-}){e_{S^1}}({{\calN_0}^+}),$$
where by ${{\calN_0}^+}$ we denote of course the positive normal
bundle to $M_0$ and by ${{\calN_0}^-}$ the negative normal bundle.
Therefore we may modify (40) by an overall (infinite!) constant
and redefine it as:
\begin{equation}
{{\mathcal{F}}_{T_{\beta}}}({\Delta})=
\int_{\tLM}{e^{{tP}/{\hbar}}}{C_{T_{\beta}}}{\Delta}=
\sum_{d\in{\check{\mathcal{K}}}}
\int_{M}{T_{\beta}}{e^{\sum_{k=1}^{r}{t_k}({{{\omega}_k}/{\hbar}}+{\int_{d}{{\omega}_k}})}}
\frac{{e_{S^1}}({{\calN_d}^+})}{{e_{S^1}}({{\calN_0}^+})}.
\end{equation}
This step may seem heretical and arbitrary to the reader (and it
certainly is)  but will be justified in the next section where an
honest mathematical proof that this "regularization" works, will
be given for the toric case.
 A similar "regularization" is used by
Givental in \cite{HomGeom}.

Now as was explained in the previous section, in the case where
the cohomology of $M$ is generated by ${H^2}(M,{\bZ})$, the
quantum $\mathcal{D}$ - module is expected to be
${\mathcal{D}}/{I_0}$ where $I_0$ is the ideal of operators
annihilating $\mathcal{F}_{T_{\beta}}({\Delta})$ for $\beta=0\dots
m$. It is more convenient to consider them all at once as follows:
Consider the following ${H^{2*}}(M,\bC)$ valued function:
\begin{equation}
F=\sum_{d\in \check{\mathcal{K}}}
{e^{\sum_{k=1}^{r}{t_k}({{{\omega}_k}/{\hbar}}+{\int_{d}{{\omega}_k}})}}
\frac{{e_{S^1}}({{\calN_d}^+})}{{e_{S^1}}({{\calN_0}^+})}
\end{equation}

Then
$$F=\sum_{\beta=1}^{m}{\mathcal{F}}_{T_{\beta}}({\Delta}){T^{\beta}}.$$

{\itshape It is important to notice that the E\"{u}ler class
$e_{S^1}({{\calN_d}^+})$ used in (41) and (42), is not well
defined, since it is an infinite divergent product} (of the
E\"{u}ler classes of the line bundles in the -possibly virtual-
splitting of ${\calN_d}^+$ according to weights of the
representation of the circle action on the fibers). What we claim
though, is that, certainly in the toric case (see next section)
and conjecturally for general symplectic manifolds {\itshape there
is a way to define the ratio of the two E\"{u}ler classes}.

Bearing this in mind, we may now formulate a conjecture for the
regularization of the ``Fourier transform'' of the Floer
fundamental cycle $\Delta$ as follows:

\begin{conjecture}
If the cohomology of $M$ is generated by classes in
${H^2}(M,{\bZ})$ then, the quantum cohomology $\mathcal{D}$ -
module of $M$ is generated by the ${H^{2*}}(M,{\bC})$ valued
function:

\begin{equation}
F=\sum_{d\in \check{\mathcal{K}}}
{e^{\sum_{k=1}^{r}{t_k}({{{\omega}_k}/{\hbar}}+{\int_{d}{{\omega}_k}})}}
\frac{{e_{S^1}}({{\calN_d}^+})}{{e_{S^1}}({{\calN_0}^+})}
\end{equation}
{\itshape or equivalently}
\begin{equation}
F={e^{({t_1}{{\omega}_1}+\dots +{t_r}{{\omega}_r})/{\hbar}}}
\sum_{d\in \check{\mathcal{K}}} {q^d}\
\frac{{e_{S^1}}({{\calN_d}^+})}{{e_{S^1}}({{\calN_0}^+})},
\end{equation}

which has components appropriate regularizations of
$\mathcal{F}_{T_{\beta}}({\Delta})$. Here ${q_k}={e^{t_k}}$,
${d_k}=\int_{d}{{\omega}_k}$ and ${q^d}={{q_1}^{d_1}}\dots
{{q_1}^{d_1}}$.

Moreover for $d=0$ we let
$\frac{{e_{S^1}}({{\calN_d}^+})}{{e_{S^1}}({{\calN_0}^+})}=1$\ .

 By ``generate the quantum cohomology $\mathcal{D}$ - module'' we mean that
the quantum cohomology $\mathcal{D}$ - module is equal to
${\mathcal{D}}/{I_0}$ where $I_0$ is the ideal of operators in
$\mathcal{D}$ that annihilate the components of $F$.

Moreover we conjecture that :
\begin{equation}
F=G,
\end{equation}

where $G$ is the function defined, in terms of the space of stable
maps, by equation (7) and by definition has the property of
generating the quantum $\mathcal{D}$ - module of $M$.

\end{conjecture}

Our next task will be to explain how to define the ratio of the
E\"{u}ler classes in the case of toric manifolds and then prove
that $F=G$ and thus indeed generates the quantum $\mathcal{D}$ -
module.

The strategy in that case will be to reformulate everything in
terms of a sequence of {\itshape finite dimensional
approximations} of the loop space, by spaces parameterizing loops
of arbitrarily large but finitely many modes (in their Fourier
expansion). We can then consider the corresponding sequence of
ratios of E\"{u}ler classes. These ratios stabilize for large
modes and we use this ``\textit{stable} \textit{ratio}'' as our
definition. We can then calculate $F$ explicitly. Finally, we
invoke the calculation of $G$ by Givental \cite{EGWI} to show that
$F=G$.

As a final observation we note that ${\mathcal{F}}_{T_{\beta}}$
for all ${\beta}=0,\dots ,m$ are elements of
${Hom_{\mathcal{D}}}({FH^*}_{S^1}({\tLM}),{\mathcal{O}})$ and
therefore define flat section of the A-connection. Denote by
${\Delta}_{\alpha}$ the Floer-Witten cycle corresponding to
$T_{\alpha}$. Then ${{\Delta}_0}=\Delta$. For a fixed $T_{\beta}$,
the functions ${\mathcal{F}}_{T_{\beta}}({{\Delta}_{\alpha}})$
simply give ${\mathcal{F}}_{T_{\beta}}$ in a basis. They should be
$${\mathcal{F}}_{T_{\beta}}({{\Delta}_{\alpha}})={s_{{\alpha},{\beta}}},$$
where $({s_{{\alpha},{\beta}}})$ is the fundamental solution
matrix (6) of the $A$ model flat section equation. In this fashion
we may identify ${\mathcal{F}}_{T_{\beta}}$ with $s_{\beta}$
namely, the flat section found in theorem (1).

\section{The Toric Case}
Our goal in this section is to formulate and prove rigorously a
version of Conjecture (1) formulated in the previous section.
First we need to describe the set up. Let $M$ be a compact,
smooth, K\"{a}hler toric variety.\footnote{Excellent references
for toric varieties are Fulton \cite{bill} for the algebraic
geometric point of view and Audin \cite {Audin} for the symplectic
side.} We choose to think of it as a symplectic quotient. To that
end, in order to define $M$ we start with an exact sequence of
lattices as in:
\begin{equation}
{0}\to {{\Bbb Z}^l}\to {{\Bbb Z}^n}\to {{\Bbb Z}^d}\to {0},
\end{equation}
where the first map is called $m$ and the second $\pi$.
Tensoring the sequence with $\bC$ and exponentiating gives a
sequence of algebraic tori :
\begin{equation}
{1}\to {{\Cstar}^l}\to {{\Cstar}^n}\to {{\Cstar}^d}\to
{1}.
\end{equation}
Now tensoring (46) with $i{\bR}$ and exponentiating gives a
sequence of real tori
\begin{equation}
{1}\to {{\bT}^l}\to {{\bT}^n}\to {{\bT}^d}\to {1}.
\end{equation}
These sequences define an embedding of ${\Cstar}^l$ into
${\Cstar}^n$ and of ${{\bT}^l}$ into ${{\bT}^n}$. Composing this
with the diagonal action of ${\Cstar}^n$ on itself defines the
action
\begin{equation}
(x_1,\dots ,x_n)\mapsto
(\prod_{j=1}^{l}{{{\lambda}_j}^{m_{j,1}}}x_1,\dots
,\prod_{j=1}^{l}{{{\lambda}_j}^{m_{j,n}}}x_n).
\end{equation}
Associated to this, there is the moment map
$${\mu}:{{\bC}^n}\to {{\bR}^l}$$ given by
\begin{equation}
{\mu}=
{\mu}_l=\check{m} \circ {{\mu}_n}={\frac{1}{2}}({\sum_{k=1}^{n}{{m_{1,k}}
{{|{x_k}|}^2}}},
\dots ,{\sum_{k=1}^{n}{{m_{l,k}}{{{|{x_k}|}^2}}}}).
\end{equation}
If $\lambda\in {{\bR}^l}$ is a regular value of $\mu$ then $M$ is
constructed by symplectic reduction as
\begin{equation}
M={M_{\lambda}}={{\mu}^{-1}}(\lambda)/{{\bT}^l}.
\end{equation}
$M$ comes equipped with the reduced symplectic form
${\omega}_{\lambda}$. For simplicity we shall just denote it by $\omega$.

$M$ is a K\"{a}hler (and at worst) orbifold. The K\"{a}hler form
is the reduction of the standard K\"{a}hler form on ${\bC}^n$.
Notice further that there is a cone in ${\bR}^l$ defined by the
conditions that it contains $\lambda$ and that the differential of
$\mu$ drops rank along its walls. Reducing at any point in the
cone gives a space topologically equivalent but with a different
K\"{a}hler form. In fact we may identify this cone with the
K\"{a}hler cone $\mathcal{K}$ of $M$.

Now recall that if $({x_1},\dots ,{x_n})$ are coordinates on
${\bC}^n$ then they can be thought of as sections of corresponding
line bundles $L_k$ over $M$ for $k=1,\dots ,n$. The divisor
$(x_k)$ is denoted by ${D_k}$. We call these divisors the toric
divisors. Let ${v_k}={\pi}(w_k)$ for $k=1,\dots n$ where
$\{{w_1},\dots ,{w_n}\}$ is the standard basis of ${\bZ}^n$, then
rational equivalences among the $D_k$ are given by the relations
\begin{equation}
\sum_{k=1}^{n}<{{e_{\nu}}^*},{v_k}>{D_k}=0\ \ for\ \  {\nu}=1,\dots
,(n-l),
\end{equation}
where $\{{{e_1}^*},\dots ,{{e_{n-l}}^*}\}$ is the dual of the
standard basis of ${\bZ}^{n-l} $. Let ${{\alpha}_k}={c_1}({L_k})$.
The ${\alpha}_k$ are Poincar\'e dual to the ${D_k}$ and (52) gives
the additive relations among them.

It is known (see eg Fulton \cite{bill})that ${H^2}(M,{\bR})$ is
spanned by the ${\alpha}_k$ and that ${H^*}(M,{\bR})$ is generated
by classes in ${H^2}(M,{\bR})$.

Now we would like to model somehow the space $\tLM$ defined in
section (3) as a covering of the space of free contractible loops
in $M$. To this end lets consider loops in ${\bC}^n$ with a
{\itshape finite but large} number of modes $2N$. To be specific
we shall consider loops which are in general of the form:
\begin{equation}
{\gamma}:{S^1}\to {{\bC}^n}\ with  \
{\gamma}({e^{i{\theta}}})=({{\gamma}_1}({e^{i{\theta}}}),\dots
,{{\gamma}_n}({e^{i{\theta}}})),
\end{equation}
where, if we let $z={e^{i{\theta}}}$, then
\begin{equation}
{{\gamma}_k}:{S^1}\to {\bC},
\end{equation}
has Fourier expansion :
\begin{equation}
{{\gamma}_k}(z)=\sum_{{\nu}=-N}^{N}{{a_{\nu}}^k}{z^{\nu}}.
\end{equation}
Our model $\tLNM$ for $\tLM$ will be defined as follows : The
space $\tLNCn$ of loops of finite modes in ${\bC}^n$ is
parametrized by the Fourier coefficients ${{a_{\nu}}^k}$ and
therefore is just ${\bC}^{n2N}$. Consider the ${\Cstar}^l$ (or
$T^l$) action on $\tLNCn$ induced by the action (49) on ${\bC}^n$
defining $M$. By this we mean that the action on all the
coefficients of ${\gamma}_k$ is the same as the action on $x_k$.
The moment map attached to this action is:
\begin{equation}
{{\mu}_N}={\frac{1}{2}}(
\sum_{{\nu}=-N}^{N}\sum_{k=1}^{n}{{m_{1,k}}{{|{{{a_{\nu}}^k}}|}^2}},
\dots
,\sum_{{\nu}=-N}^{N}\sum_{k=1}^{n}{{m_{l,k}}{{|{{{a_{\nu}}^k}}|}^2}}).
\end{equation}

Define $\tLNM$ as
\begin{equation}
\tLNM={{\mu}_N}^{-1}({\lambda})/{T}^l,
\end{equation}
where we have identified the K\"{a}hler cones of $M$ and $\tLNM$. We can
do this since the subsets of ${\bR}^l$ where $\mu$ and ${\mu}_N$ drop
rank are clearly the same. This is just because $\mu$ drops rank at
some value if and only if some
 homogeneous coordinates are forced to be zero.
At the same value ${\mu}_N$
 drops rank since the corresponding sums of squares  are forced to be
zero which in turn forces each of the squares to be zero.

Now in general the K\"{a}hler cone $\mathcal{K}$  will not necessarily
be simplicial, but
it can of course be subdivided, to simplicial cones. Pick such a subdivision
and consider the simplicial cone containing the value
$\lambda$.

Let $\{{\omega_1},\dots {\omega_l}\}$ be the basis of that cone such
that
\begin{equation}
{\omega}=\sum_{j=1}^{l}{{\lambda}_j}{{\omega}_j}.
\end{equation}

The fact that the cone is simplicial means that
$\{{\omega_1},\dots {\omega_l}\}$ is a basis of ${H^2}(M,{\bZ})$.

Moreover we have that

\begin{equation}
{{\alpha}_k}=\sum_{j=1}^{l}{m_{j,k}}{{\omega}_j}.
\end{equation}

If $d$ is an element of $H_2(M,{\bZ})$ then we let
\begin{equation}
{d_j}=\int_{d}{\omega_j}.
\end{equation}
In that case we may identify $d$ with the vector $({d_1},\dots
,{d_l})$.

Next we need to consider the action functional $H_N$ associated to
$\tLNM$. Recall that the action functional assigns to a pair
(loop, contracting disc) the symplectic area of the contracting
disc. Recall also that we have fixed the standard K\"{a}hler form
on ${\bC}^n$ which is
\begin{equation}
{{\omega}_0}=\sum_{k=1}^n {{ds}_{k}}\wedge {{dt}_{k}} =
\frac{i}{2}\ \sum_{k=1}^n {{dx}_k}\wedge {d}\overline{{x}_k},
\end{equation}
where ${x_k}={s_k}+i{t_k}$.
 Define first
\begin{equation}
{H_N}({\gamma})=\frac{1}{2\pi}\int_{{\gamma}(S^1)}\sum_{k=1}^{n}{s_k}{dt_k}=
\frac{1}{2\pi}\sum_{k=1}^{n}\int_{S^1}{{\gamma}^*}({s_k}{dt_k})=
\frac{1}{2\pi}\sum_{k=1}^{n}\int_{D}{{u}^*}(\frac{i}{2}{{dx}_k}\wedge
{d}\overline{{x}_k}),
\end{equation}
where $\gamma$ is given by (53). In other words $H_N$ is the
(normalized) action functional for loops in ${\bC}^n$ or rather
$\tLNCn$ to be exact.

An elementary calculation shows that
\begin{equation}
{H_N}({\gamma})=\frac{1}{2}\sum_{{\nu}=-N}^{N}{\nu}({{|{{{a_{\nu}}^1}}|}^2}+\dots
+{{|{{{a_{\nu}}^n}}|}^2}).
\end{equation}
To see this, it's enough to notice that if a loop $\gamma:
 {S^1} \to\bC$ is given by $\gamma(e^{i \theta})=e^{i k \theta}$
 then $u:D\to \bC$ such that $u(re^{i \theta})={r^k}e^{\imath k
\theta}$ contracts that loop. Moreover if $x=re^{i \theta}$ is a
coordinate on $\bC$ then
$$\int_Du^*(dx\wedge d\overline{x})=\int_D(du\wedge d\overline{u})=
\int_0^{2\pi}\int_0^1{-2}i {k^2} r^{2k-1}dr\wedge d \theta= -2\pi
i k$$
 Now $H_N$ is thus far defined on ${\bC}^{n2N}$, but since it is
invariant under the $T^l$ action, it actually drops to a function
on $\tLNM$. We will still call that function by the same name
$H_N$, and it is our action functional.

Consider next the $S^1$ action on $\tLNM$. It is induced by rotation
on the source circle, namely by the action
${e^{i{\theta}}}\mapsto {e^{i({\theta}+{\phi})}}$ .
This action induces an action on the Fourier coefficients of a loop
$\gamma$ by
\begin{equation}
{{{a_{\nu}}^k}}\mapsto {e^{ik{\phi}}}{{{a_{\nu}}^k}}.
\end{equation}
It's clear that $H_N$ is the Hamiltonian function corresponding to
this action on $\tLNCn$ and consequently
 on $\tLNM$. This
is in accordance with the general theory. As we saw in equation
(16), the action functional is indeed the Hamiltonian of the
circle action.

As a consequence the fixed components of the circle action on
$\tLNM$ coincide with the critical manifolds of $H_N$. We expect
those to be copies of $M$ and to correspond to homology classes
$d\in {H_2}(M,{\bZ})$. Recall that we may identify the class $d$
with its period vector $({d_1},\dots ,{d_l})$ as in (60). Now the
fixed components of the circle action can be identified as
follows: The action of $T^l$ on
${{\mu}^{-1}}({\lambda})\in\tLNCn={{\bC}^{n2N}}$ that defines
$\tLNM$ is induced by the action in (49). When we take
${{\lambda}_1}={z^{d_1}},\dots ,{{\lambda}_l}={z^{d_l}}$, this
becomes an $S^1$ action. Components in ${{\mu}^{-1}}({\lambda})$
where the $S^1$ action from (64) coincides with the one appearing
as a one parameter subgroup of the $T^l$ action as above, will
lead to fixed components in $\tLNM={{\mu}^{-1}}({\lambda})/{T^l}$.
Thus making the substitution
$${{\lambda}_1}={z^{d_1}},\dots ,{{\lambda}_l}={z^{d_l}}$$
in (49) shows immediately that loops of the form
\begin{equation}
{\gamma}(z)=
({{{a^1}_{\sum_{j=1}^{l}{m_{j,1}}{d_j}}}}{z^{\sum_{j=1}^{l}{m_{j,1}}{d_j}}},\dots ,
{{{a^n}_{\sum_{j=1}^{l}{m_{j,n}}{d_j}}}}{z^{\sum_{j=1}^{l}{m_{j,n}}{d_j}}}),
\end{equation}
form a fixed component of the circle action on $\tLNM$. The $T^l$
action on ${\bC}^{n2N}$ restricted to loops in ${\bC}^n$ of the
form (65), restricts to the action (49) defining $M$. Therefore
the reduction of the space of loops of the form (65) will indeed
be exactly a copy of $M$. We shall name this component $M_d$.
Notice that there is a more illuminating way to write (65).
According to (59) and (60)  we have
\begin{equation}
{\sum_{j=1}^{l}{m_{j,k}}{d_j}}=\int_{d}{{\alpha}_k},
\end{equation}
where ${{\alpha}_k}$, as we said earlier, is Poincar\'e dual to the toric
divisor $D_k$.
Therefore we see that $M_d$ consists of loops of the form :
\begin{equation}
{\gamma}(z)=({{{a^1}_{{\int_{d}{{\alpha}}}}}}{z^{{\int_{d}{{\alpha}_1}}}},\dots
,{{{a^n}_{{\int_{d}{{\alpha}_1}}}}}{z^{{\int_{d}{{\alpha}_1}}}}).
\end{equation}
This immediately tells us that in order to be able to study $M_d$
we must take $N\geq
max\{\int_{d}\alpha_1,\dots,\int_{d}\alpha_n\}$ which we will
assume from now on whenever discussing $M_d$. To recapitulate the
set up, up to now we have defined the spaces $\tLNM$ which as
$N\to \infty$ approximate $\tLM$ and action functionals
\begin{equation}
{H_N}:{\tLNM}\to {\bR}.
\end{equation}
We have also described the critical manifolds $M_d$ of $H_N$,
which are copies of $M$ and of course coincide with the fixed
components of the circle action on $\tLNM$.

Moreover notice that since $H_N$ is the Hamiltonian of an $S^1$
action it follows from general theory that it is a perfect
Morse-Bott function. This is explained for example in Audin
\cite{Audin}. In our case it also obvious from (63) which shows
that indices of $H_N$ are even numbers. The so called lacunary
principle (see e.g. Bott \cite{Bott}) then guarantees that $H_N$
is perfect. It's also clear that $H_N$ is non-degenerate in the
normal directions.

 Now let ${\calN}_{d,N}$
denote that normal bundle to ${M_d}\subset {\tLNM}$. Let
${{\calN}^+}_{d,N}$ and ${{\calN}^-}_{d,N}$ denote the positive
and negative normal bundles of ${M_d}\subset {\tLNM}$. Obviously
we have ${{\calN}_{d,N}}={{{\calN}^-}_{d,N}}\bigoplus
{{{\calN}^+}_{d,N}} $. We shall first prove the following
proposition.
\begin{proposition}
Let $M$ be toric manifold. For every class $d\in {H_2}(M,{\bZ})$,
there is an integer $N(d)$ such that the ratio of equivariant
E\"{u}ler classes,
\begin{equation}
\frac{{e_{S^1}}({{{\calN}^+}_{d,N}})}{{e_{S^1}}({{{\calN}^+}_{0,N}})},
\end{equation}
remains constant for all $N\geq N(d)$.
\end{proposition}

In other words this ratio of E\"{u}ler classes stabilizes. This
allows us to define the {\itshape stable ratio}.
\begin{definition}
Define the stable ratio
$\frac{{e_{S^1}}({{{\calN}^+}_d})}{{e_{S^1}}({{{\calN}^+}_0})}$ to be
the common ratio
$\frac{{e_{S^1}}({{{\calN}^+}_{d,N}})}{{e_{S^1}}({{{\calN}^+}_{0,N}})}$
for all $N\geq N(d)$
\end{definition}

After we have this and keeping the notation from above we will
prove the following version of Conjecture 1.
\begin{theorem}
Let $M$ be a smooth toric variety of Picard number $l$, and assume
that for every $d \in {\check{\mK}}-\{0\}$ we have $\int_d{c_1}({T_M})> 0$.
Let
\begin{equation}
F={e^{({t_1}{{\omega}_1}+\dots +{t_l}{{\omega}_l})/{\hbar}}}
\sum_{d\in \check{\mathcal{K}}}
{q^d}\ \frac{{e_{S^1}}({{{\calN}^+}_d})}{{e_{S^1}}({{{\calN}^+}_0})},
\end{equation}
then
\begin{equation}
F=G,
\end{equation}
where $$G={e^{({t_1}{{\omega}_1}+\dots
+{t_l}{{\omega}_l})/{\hbar}}}(1+\sum_{d\in
{\check{\mathcal{K}}},d\ne
0}{q^d}{{ev_1}_*}(\frac{1}{{\hbar}-c}))$$ is the function defined
in (7) (since ${p_i}={\omega_i}$ and ${q_i}=e^{t_i})$. Here
${ev_1}:{M_{2,d}}\to M$ is evaluation at the first marked point
and $c$ is the Chern class of the line bundle over $M_{2,d}$ which
is the universal cotangent line at the second marked point.
Moreover, since $G$ is the function that generates the Quantum
cohomology $\mathcal{D}$ - module of $M$, it follows that $F$ does
too.
\end{theorem}

We begin with the proof of proposition (3). We need to analyze in
detail the normal bundle ${{\calN}_{d,N}}$. One way to do this,
would be to analyze the normal bundle by considering the affine
patches that cover $\tLNM$ coming from toric geometry. This is
feasible but not so convenient so we take a different approach.
Recall that $\tLNM={{{\mu}_N}^{-1}}({\lambda})/{T^l}$, where
$\lambda \in \bR^l$ has components ${\lambda}=({{\lambda}_1},\dots
,{{\lambda}_l})$. Therefore according to (56) we  have
\begin{equation}
\frac{1}{2}\sum_{{\nu}=-N}^{N}\sum_{k=1}^{n}{{m_{1,k}}{{|{{{a_{\nu}}^k}}|}^2}}={{\lambda}_1},
\dots ,
\frac{1}{2}\sum_{{\nu}=-N}^{N}\sum_{k=1}^{n}{{m_{l,k}{{|{{{a_{\nu}}^k}}|}^2}}}=
{{\lambda}_l}.
\end{equation}
Now the action functional $H_N$ can be restricted to
${{{\mu}_N}^{-1}}({\lambda})$ by using the relations (72). Once we
do this, then we have the function on $\tLNM$ since $H_N$ is
invariant under the $T_l$ action. In order to look in the normal
directions of $M_d$ and using the description of $M_d$ found in
(67),  we work as follows:
\begin{equation}
2{H_N}({\gamma})=\sum_{{\nu}=-N}^{N}{\nu}({{|{{{a_{\nu}}^1}}|}^2}+\dots
+{{|{{{a_{\nu}}^n}}|}^2})=
\end{equation}
$$=\sum_{{\nu}\ne {\int_{d}{{\alpha}_{1}}}}{\nu}{{|{{{a_{\nu}}^1}}|}^2}+\dots
+\sum_{{\nu}\ne
{\int_{d}{{\alpha}_{n}}}}{\nu}{{|{{{a_{\nu}}^n}}|}^2}+
(\int_{d}{{\alpha}_1}){{|{{a^1}_{\int_{d}{{\alpha}_1}}}|}^2}+
\dots +
(\int_{d}{{\alpha}_n}){{|{{a^n}_{\int_{d}{{\alpha}_n}}}|}^2}.$$
Now substituting  (66) in (73) we find :
\begin{equation}
2{H_N}=
\sum_{{\nu}\ne {\int_{d}{{\alpha}_{1}}}}{\nu}{{|{{{a_{\nu}}^1}}|}^2}+\dots
+\sum_{{\nu}\ne {\int_{d}{{\alpha}_{n}}}}{\nu}{{|{{{a_{\nu}}^n}}|}^2}+
\sum_{j=1}^{l}{m_{j,1}}{d_j}{{|{{a^1}_{\int_{d}{{\alpha}_1}}}|}^2}+ \dots +
\sum_{j=1}^{l}{m_{j,n}}{d_j}{{|{{a^n}_{\int_{d}{{\alpha}_n}}}|}^2}.
\end{equation}
Rearranging this sum gives :
\begin{equation}
2{H_N}=
\sum_{{\nu}\ne {\int_{d}{{\alpha}_{1}}}}{\nu}{{|{{{a_{\nu}}^1}}|}^2}+\dots
+\sum_{{\nu}\ne {\int_{d}{{\alpha}_{n}}}}{\nu}{{|{{{a_{\nu}}^n}}|}^2}+
{d_1}\sum_{k=1}^{n}{m_{1,k}}{{|{{a^k}_{\int_{d}{{\alpha}_k}}}|}^2}+ \dots +
{d_l}\sum_{k=1}^{n}{m_{l,k}}{{|{{a^k}_{\int_{d}{{\alpha}_k}}}|}^2}.
\end{equation}
Now we may use (72) to obtain

$$2{H_N}=2\sum_{j=1}^{l}{d_j}{{\lambda}_j}+
\sum_{{\nu}\ne {\int_{d}{{\alpha}_{1}}}}{\nu}{{|{{{a_{\nu}}^1}}|}^2}+\dots
+\sum_{{\nu}\ne {\int_{d}{{\alpha}_{n}}}}{\nu}{{|{{{a_{\nu}}^n}}|}^2}-$$
$$-\sum_{k=1}^{n}\sum_{{\nu}\ne {\int_{d}{{\alpha}_k}}}{d_1}{m_{1,k}}
{{|{{a_{\nu}}^k}|}^2}-\dots - \sum_{k=1}^{n}\sum_{{\nu}\ne
{\int_{d}{{\alpha}_k}}}{d_l}{m_{l,k}} {{|{{a_{\nu}}^k}|}^2}.$$
Rearranging the sum once more we find :
$${H_N}=\sum_{j=1}^{l}{d_j}{{\lambda}_j}+
\frac{1}{2}\sum_{{\nu}\ne
{\int_{d}{{\alpha}_{1}}}}({\nu}-\sum_{j=1}^{l}{d_j}{m_{j,1}})
{{|{{{a_{\nu}}^1}}|}^2}+\dots +
\frac{1}{2}\sum_{{\nu}\ne
{\int_{d}{{\alpha}_{n}}}}({\nu}-\sum_{j=1}^{l}{d_j}{m_{j,n}})
{{|{{{a_{\nu}}^n}}|}^2}$$
and finally
\begin{equation}
{H_N}=\int_{d}{\omega}+
\frac{1}{2}\sum_{{\nu}\ne {\int_{d}{{\alpha}_{1}}}}({\nu}-\int_{d}{{\alpha}_1})
{{|{{{a_{\nu}}^1}}|}^2}+\dots +
\frac{1}{2}\sum_{{\nu}\ne {\int_{d}{{\alpha}_{n}}}}({\nu}-\int_{d}{{\alpha}_n})
{{|{{{a_{\nu}}^n}}|}^2}.
\end{equation}

Notice that the fact that
$\sum_{j=1}^{l}{d_j}{{\lambda}_j}=\int_{d}{\omega},$ follows from
(58) and (60). Moreover since  $H_N$ is quadratic, computing the
Hessian is immediate. Recall also from (67) that
\begin{equation}
{M_d}= ({{a^1}_{{\int_{d}{{\alpha}_1}}}},\dots
,{{a^n}_{{\int_{d}{{\alpha}_1}}}})/{T^l},
\end{equation}
where it is implied that all the other coordinates are zero. The
normal bundle of the subset of ${\bC}^{n2N}$ whose quotient by the
$T^l$ action is $M_d$, is trivial of course and has fiber
coordinates given by all the variables in ${\bC}^{n2N}$ except for
$({a_{{\int_{d}{{\alpha}_1}}}},\dots
,{a_{{\int_{d}{{\alpha}_1}}}})$. The normal bundle of $M_d$ is the
quotient of the normal bundle to the subset
$({a_{{\int_{d}{{\alpha}_1}}}},\dots
,{a_{{\int_{d}{{\alpha}_1}}}})$. Therefore if we let
${\tilde{L}}_{k,{\nu}}$ be the bundle over $\tLNM$ associated to
${a_{\nu}}^k$ and define
\begin{equation}
{E_d}=\bigoplus_{{\nu}\ne {\int_{d}{{\alpha}_{1}}}}{{\tilde{L}}_{1,{\nu}}}\oplus \dots
\oplus \bigoplus_{{\nu}\ne {\int_{d}{{\alpha}_{n}}}}{{\tilde{L}}_{n,{\nu}}},
\end{equation}
then
\begin{equation}
{{\mathcal{N}}_{d,N}}={E_d}|_{M_d}.
\end{equation}
Of course ${\tilde{L}}_{k,{{\nu}_1}}$ is isomorphic to
${\tilde{L}}_{k,{{\nu}_2}}$ for all ${\nu}_1$ and ${\nu}_2$ if we
don't take into account the $S^1$ action but they are different as
$S^1$ equivariant bundles according to (64). Now let
\begin{equation}
{L_{k,{\nu}}}={{\tilde{L}}_{k,{\nu}}}|_{M_d},
\end{equation}
then we have that
\begin{equation}
{{L_{k,{\nu}}}\ {\rm isomorphic\ {to}}\ {L_k}} \ {\rm for\ all}\
{\nu},
\end{equation}
where again the isomorphism is taken in the usual sense, not the $S^1$
equivariant.

As a consequence of (79) we have that
\begin{equation}
{{\mathcal{N}}_{d,N}}=
\bigoplus_{{\nu}\ne {\int_{d}{{\alpha}_{1}}}}{L_{1,{\nu}}}\oplus \dots
\oplus \bigoplus_{{\nu}\ne {\int_{d}{{\alpha}_{n}}}}{L_{n,{\nu}}}.
\end{equation}

It is now straight forward to compute the $S^1$ equivariant
E\"{u}ler class of each of the line bundles. The key point for us,
is to read off from the expression of the action functional in
(76) the {\itshape positive} normal bundle. Clearly we have that
if
\begin{equation}
{{E^+}_d}=\bigoplus_{{\nu}> {\int_{d}{{\alpha}_{1}}}}{{\tilde{L}}_{1,{\nu}}}\oplus \dots
\oplus \bigoplus_{{\nu}> {\int_{d}{{\alpha}_{n}}}}{{\tilde{L}}_{n,{\nu}}},
\end{equation}
then
\begin{equation}
{{{\mathcal{N}}^+}_{d,N}}={{{E^+}_d}|{_{M_d}}}=
\bigoplus_{{\nu}> {\int_{d}{{\alpha}_{1}}}}^N{L_{1,{\nu}}}\oplus \dots
\oplus \bigoplus_{{\nu}> {\int_{d}{{\alpha}_{n}}}}^N{L_{n,{\nu}}}.
\end{equation}

Finally since the $S^1$ action is given by (64) and using (80) we
find that
\begin{equation}
{e_{S^1}}({L_{k,{\nu}}})={c_1}({L_{k,{\nu}}})+{\nu}{\hbar}={{\alpha}_k}+{\nu}{\hbar},
\end{equation}
and therefore
\begin{equation}
{e_{S^1}}({{{\mathcal{N}}^+}_{d,N}})=
{\prod_{{\nu}>
{\int_{d}{{\alpha}_{1}}}}^N({{\alpha}_1}+{\nu}{\hbar})}\dots
{\prod_{{\nu}>
{\int_{d}{{\alpha}_{n}}}}^N({{\alpha}_n}+{\nu}{\hbar})}.
\end{equation}
As a special case it follows that
\begin{equation}
{e_{S^1}}({{{\mathcal{N}}^+}_{0,N}})=
{\prod_{{\nu}>0}^N({{\alpha}_1}+{\nu}{\hbar})}\dots
{\prod_{{\nu}>0}^N({{\alpha}_n}+{\nu}{\hbar})}.
\end{equation}
We may now compute the ratio
\begin{equation}
\frac{{e_{S^1}}({{{\mathcal{N}}^+}_{d,N}})}
{{e_{S^1}}({{{\mathcal{N}}^+}_{0,N}})}=
\frac{1}
{{\prod_{{\nu}=1}^{{\int_{d}{{\alpha}_{1}}}}({{\alpha}_1}+{\nu}{\hbar})}
\dots
{\prod_{{\nu}=1}^{{\int_{d}{{\alpha}_{n}}}}({{\alpha}_n}+{\nu}{\hbar})}}.
\end{equation}
It follows that indeed the ratio
$\frac{{e_{S^1}}({{{\mathcal{N}}^+}_{d,N}})}
{{e_{S^1}}({{{\mathcal{N}}^+}_{0,N}})}$ is independent of $N$ as
long as $N$ is greater than
$N(d)=max\{{{\int_{d}{{\alpha}_{1}}}},\dots ,
{{\int_{d}{{\alpha}_{n}}}}\}$. This concludes the proof of
Proposition (3).

The proof of Theorem (4) follows now easily from the proof of
Proposition (3). First note that the \textit{stable}
\textit{ratio} is
\begin{equation}
\frac{{e_{S^1}}({{{\mathcal{N}}^+}_d})}
{{e_{S^1}}({{{\mathcal{N}}^+}_0})}= \frac{1}
{{\prod_{{\nu}=1}^{{\int_{d}{{\alpha}_{1}}}}({{\alpha}_1}+{\nu}{\hbar})}
\dots
{\prod_{{\nu}=1}^{{\int_{d}{{\alpha}_{n}}}}({{\alpha}_n}+{\nu}{\hbar})}}.
\end{equation}

We may now compute the function $F$ of (70). We find that
\begin{equation}
F={e^{({t_1}{{\omega}_1}+\dots +{t_l}{{\omega}_l})/{\hbar}}}
\sum_{d\in \check{\mathcal{K}}}
{q^d}\ \frac{{e_{S^1}}({{{\calN}^+}_d})}{{e_{S^1}}({{{\calN}^+}_0})}=
\end{equation}
$$={e^{({t_1}{{\omega}_1}+\dots +{t_l}{{\omega}_l})/{\hbar}}}
\sum_{d\in \check{\mathcal{K}}}
{q^d}\ \frac{1}
{{\prod_{{\nu}=1}^{{\int_{d}{{\alpha}_{1}}}}({{\alpha}_1}+{\nu}{\hbar})}
\dots
{\prod_{{\nu}=1}^{{\int_{d}{{\alpha}_{n}}}}({{\alpha}_n}+{\nu}{\hbar})}},
$$
where ${q^d}$ stands for
${{q_1}^{d_1}}\dots {{q_l}^{d_l}}={e^{{t_1}{d_1}+\dots +{t_l}{d_l}}}$ as usual.

Finally, according to Givental's computation of the function
$$G={e^{plnq/{\hbar}}}(1+\sum_{d\in {\check{\mathcal{K}}},d\ne
0}{q^d}{{ev_1}_*}(\frac{1}{{\hbar}-c}))$$ in \cite{toric} (Theorem
(0.1), page (3) and its corollary: Example (a) page (4)) we have
that, if $\int_{d}{c_1}(TM)>0$ for all $d \in {\check{\mK}}$ and
$d\neq 0$, then  the function $F$ as computed above is indeed
equal to the function $G$ and therefore it generates the Quantum
$\mathcal{D}$ - module. This concludes the proof of theorem (4).
It may finally be useful to consider a simple example in order to
clarify things a bit more.

\begin{example}
Let us consider the simplest example which is the complex
projective space ${\bP}^n$. Let $\omega$ be the class dual to a
hyperplane. The K\"{a}hler cone is a half line and is generated by
$\omega$. The toric divisor classes ${\alpha}_i$ are also all
equal to the class dual to a hyperplane. For $d$ in
${H_2}({{\bP}^n},{\bZ})$ let
${d_1}={\int_{\omega}d}=\int_{d}{{\alpha}_i}$. Let
${q_1}={e^{t_1}}$. Then the function $F$ of (70) becomes

$$F={e^{{{t_1}{\omega}}/{\hbar}}}
\sum_{d\in \check{\mathcal{K}}}
{q^d}\ \frac{{e_{S^1}}({{{\calN}^+}_d})}{{e_{S^1}}({{{\calN}^+}_0})}
={e^{{t_1}{\omega}/{\hbar}}}
\sum_{{d_1}=0}^{\infty}
{{q_1}^{d_1}}\ \frac{1}
{\prod_{{\nu}=1}^{d_1}{{({\omega}+{\nu}{\hbar})}^{n+1}}}.
$$

We may now expand $F$ in the basis
$\{1,{\omega},{{\omega}^2},\dots ,{{\omega}^n}\}$ :

$$F=\sum_{i=0}^{n}{f_i}({e^{t_1}},{\hbar}){{\omega}^i}.$$

Let $<a,b>=\int_{{\bP}^n}a\wedge b$ where $a$ and $b$ are cohomology
classes in ${\bP}^n$.
Clearly we have ${f_i}=<F,{{\omega}^{n-i}}>$
Moreover notice that
$<a,b>=Res_0 ab\frac{d\omega}{{\omega}^{n+1}}$.
Therefore
$${f_i}={Res_0}
\sum_{{d_1}=0}^{\infty}
{{q_1}^{d_1}}\ \frac{{{\omega}^{n-i}}{e^{{t_1}{\omega}/{\hbar}}}}
{\prod_{{\nu}=1}^{d_1}{{({\omega}+{\nu}{\hbar})}^{n+1}}}
\frac{d\omega}{{\omega}^{n+1}}.$$

The easiest one to compute is $f_0$ :
$${f_0}({e^{t_1}},\hbar)={Res_0}
\sum_{{d_1}=0}^{\infty} {{q_1}^{d_1}}\
\frac{{{\omega}^{n}}{e^{{t_1}{\omega}/{\hbar}}}}
{\prod_{{\nu}=1}^{d_1}{{({\omega}+{\nu}{\hbar})}^{n+1}}}
\frac{d\omega}{{\omega}^{n+1}}.$$

Therefore
$${f_0}({e^{t_1}},\hbar)=
\sum_{{d_1}=0}^{\infty}
\frac{{e^{{{t_1}{d_1}}}}}{{{\hbar}^{{d_1}(n+1)}}{({{d_1}!})^{n+1}}}.$$
The function ${f_0}$ is annihilated by the differential operator
$R({\hbar}\frac{\p}{{\p}{t_1}},{e^{t_1}},{\hbar})=
{({{\hbar}\frac{\p}{{\p}{t_1}}})^{n+1}}- {e^{t_1}}$. The quantum
$\mathcal{D}$ - module of ${\bP}^n$ is the Heisenberg algebra
modulo the ideal generated by $R$. Finally the corresponding
relation in the quantum ring of ${\bP}^n$ is $R(p,q,0)=0$ i.e.,
${p^{n+1}}=q$. Indeed, the quantum cohomology of ${\bP}^n$ is
${\bC}[p,q]/({p^{n+1}}=q)$ where $p^{n+1}$ is computed by the
quantum multiplication and $p$ is the class of the hyperplane. For
a computation of the quantum cohomology in terms of the space of
stable maps see for example \cite{EGWI}. Notice also that the rest
of the $f_i$ are also annihilated by $R$. In fact we get a
complete basis of solutions of the equation $R=0$.

It is always true that the components of $F$ satisfy the same
differential equation since they all come from the Fourier transform
of the same cycle, i.e., the Floer semi-infinite cycle $\Delta$.

\end{example}
\section{Comments and further problems}

Let us note first, that Jones and Petrack in \cite{JP} have
actually constructed an extension of equivariant cohomology that
works well in infinite dimensional settings and in particular for
the $S^1$ action on the loop space. It would be very interesting
to try and use it in this case.

Another observation follows directly from looking at the formulas
for $F$ (70) and $G$ (7). It is clear that at least in the toric
case
\begin{equation}
\frac{{e_{S^1}}({{{\calN}^+}_d})}{{e_{S^1}}({{{\calN}^+}_0})}=
{{ev_1}_*}(\frac{1}{{\hbar}-c}).
\end{equation}
This kind of formula needs to be understood in more general
symplectic manifolds and will possibly lead to a way to regularize
the ratio of E\"{u}ler classes in general. We have made some
progress in this direction and plan to report it in a subsequent
paper.

Finally we should mention that another interesting point of view
is that of Givental in \cite{HomGeom}. There, he thinks of the
integrals (39) as integrals over the Poincar\'{e} dual of
$\Delta$. As has already been observed this is geometrically the
cycle of all loops which are boundaries of holomorphic discs. Now
we have also observed that the Cauchy-Riemann equation is the flow
equation (18) for the action functional. The equation remains
unchanged if we change $H$ with $-H$ and $J$ with $-J$. Now the
change of the sign of $H$, in the $S^1$ equivariant theory, is
equivalent to changing the sign of $\hbar$. Moreover the change in
the sign of $J$ is equivalent to changing the symplectic form
$\omega$ to $-\omega$. Notice finally that Morse theory with $-H$
switches accenting with descending cells and the intersection of
accenting and descending cells gives the spaces of flow lines
connecting fixed manifolds. In our case these are spaces of
holomorphic spheres. This kind of reasoning can be beautifully
incarnated in concrete formulas. For example in the case of
$\bP^n$, Givental \cite{HomGeom} proves that:
$$\sum_{d=0}^{\infty}{e^{d\tau}}\int_{M_d}e^{(t-\tau)(H+\omega/\hbar)}=
\hbar^{1-n}\langle F(t,\hbar),\sigma F(\tau,-\hbar)\rangle,$$
where $\langle,\rangle$ denotes the intersection pairing in the
cohomology ring $\bC[P]/(P^n)$ of $\bP^n$ and $\sigma$ is the
automorphism $\sigma(P)=-P$. Here $M_d$ is the projective space of
degree $d$ polynomial maps of $\bP^1$ to $\bP^n$, $\omega$ is the
standard Fubini K\"{a}hler form on $M_d$ and $H$ is the
Hamiltonian corresponding to the $S^1$ action induced on $M_d$ by
rotation of $\bP^1$. Similar theorems exist for the function $G$
when the space of stable maps is used (see Givental
\cite{EGWI},\cite{toric}).

\bibliographystyle{amsplain}
\bibliography{loops}

\end{document}